\RequirePackage{fix-cm}

\documentclass[smallextended]{svjour3}       
\smartqed  

\usepackage{graphicx}
\usepackage{setspace}

\usepackage{float}
\usepackage{breakcites}
\usepackage{epsfig}
\usepackage{amssymb}
\usepackage{amsmath}
\usepackage{hyperref}
\usepackage{verbatim}
\usepackage{subcaption}
\usepackage{tikz}
\usepackage{color}

\usetikzlibrary{arrows}
\usepackage[bottom=0.95in,top=0.95 in]{geometry}
\newcommand\given[1][]{\:#1\vert\:}
\DeclareMathOperator{\E}{\mathbb{E}}
\tikzset{cross/.style={cross out, draw=black, minimum size=2*(#1-\pgflinewidth), inner sep=0pt, outer sep=0pt},
cross/.default={5pt}}
\usetikzlibrary{shapes.misc}
\newcommand{\be}{\begin{eqnarray}}
\newcommand{\ee}{\end{eqnarray}}
\newcommand{\bea}{\begin{eqnarray*}}
\newcommand{\eea}{\end{eqnarray*}}

\newcommand{\s}{\mathbb{S}}
\usepackage{nicefrac}
\numberwithin{equation}{section}

\usepackage{pst-node,pst-plot}
\pstVerb{realtime srand}
\psset{plotpoints=50}

\begin{document}

\title{Approximations for the  boundary crossing probabilities of moving sums of random variables}

\titlerunning{Approximations of boundary crossing probabilities}

\author{Jack Noonan    \and Anatoly Zhigljavsky 
}

\institute{Jack Noonan  \at
              School of Mathematics, Cardiff University, Cardiff, CF24 4AG, UK \\
              \email{Noonanj1@cf.ac.uk}           
           \and
           Anatoly Zhigljavsky \at
              School of Mathematics, Cardiff University, Cardiff, CF24 4AG, UK\\
              \email{ZhigljavskyAA@cardiff.ac.uk}
}

\date{Received: date / Accepted: date}

\maketitle

\begin{abstract}
In this paper we study approximations for the  boundary crossing probabilities of moving sums of i.i.d. normal random variables. We approximate a discrete time problem with a continuous time problem allowing us to apply established theory for stationary Gaussian processes. By then subsequently correcting approximations for discrete time, we show that the developed approximations are very accurate even for a small window length. Also, they have high accuracy when the original r.v. are not exactly normal and when the weights in the moving window are not all equal. We then provide accurate and simple approximations for ARL, the average run length until crossing the boundary.


\end{abstract}

\section{Introduction: Statement of the problem}
\label{sec:prob-state}

Let $\varepsilon_1,\varepsilon_2,\ldots$ be a sequence of i.i.d. normal random variables (r.v.) with mean $\theta$ and variance $\sigma^2>0$.
For a fixed positive integer~$L$, the moving sums are defined by
\begin{eqnarray}
S_{n,L}:= \sum_{j=n+1}^{n+L} \varepsilon_j\, \;\; (n=0,1, \ldots) \label{eq:sumsq2}.
\end{eqnarray}
The sequence of the moving sums \eqref{eq:sumsq2} will be denoted by $\s$ so that $\;\s=\{ S_{0,L}, S_{1,L}, \ldots \}$.

The main aim of this paper is development of accurate  approximations
 for the boundary crossing probability (BCP) for
the maximum of the moving sums:
\be
{\mathcal{P}  }_{\s}(M,H,L) := {\rm Pr}\left(\max_{n=0,1,\ldots,M} S_{n,L} \geq H\right),
 \label{eq:prob-S}
\ee
where $M$ is a given positive integer and $H$ is a fixed  threshold. Note that the total number of r.v. $\varepsilon_i$ used in~\eqref{eq:prob-S} is $M+L$ and ${\cal P}_{\s}(M,H,L) \to 1$ as $M \to \infty$, for all $H$ and $L$.
We will mostly be interested in deriving accurate approximations when $M\geq L$. The case of $M \leq L$ is much simpler and is comprehensively covered in {  \cite[Section 3]{AandZ2019},} see  Section~\ref{T1T2} for discussion.

Developing accurate approximations for the BCP ${\cal P}_{\s}(M,H,L) $  for generic parameters $H$, $M$ and $L$  is very important in various areas of statistics,  predominantly in applications related to change-point detection; see, for example, papers  \cite{Bau2,Chu,Glaz2012,MZ2003,Xia} and especially books \cite{glaz2001scan,glaz2009scan2}. Engineering applications of MOSUM (moving sums charts) are extremely important and have been widely discussed in literature; see e.g. \cite{Chu,eiauer1978use,glaz2001scan,glaz2009scan2,waldmann1986bounds}.
     The BCP ${\cal P}_{\s}(M,H,L) $ is an ($M+1$)-dimensional integral and therefore direct evaluation of this BCP is hardly possible even with modern software.

To derive approximations for the  BCP \eqref{eq:prob-S} one can use standard tools and approximate the sequence of moving sums  with a continuous-time process and then use some continuous-time approximations, see e.g. \cite{Haiman}; these approximations, however, are not accurate especially for small window length $L$; see discussion in Section~\ref{CSA_sim_study}.
There is, therefore, a need for derivation of specific approximations for the BCP \eqref{eq:prob-S}. Such a need was well
understood in the statistical community and indeed very accurate approximations for the BCP  and the Average Run Length (ARL) have been developed in a series of quality papers
by J. Glaz and coauthors, see for example \cite{Glaz_old,Glaz2012, wang2014variable,wang2014multiple} (the methodology  was  also  extended
to the case when $\varepsilon_j$ are integer-valued r.v., see  \cite{glaz1991tight}).
 We will call these approximations `Glaz approximations' by the name of the main author of these papers; they will be formally written down in Sections~\ref{Glaz_ref} and \ref{ARL_section}.

 The accuracy of the approximations developed in the present paper is very high and similar  to the Glaz approximations; this is discussed in Sections~\ref{sim_study2} and \ref{ARL_section}. The methodologies of derivation of  Glaz approximations and  the approximations of this paper are very different.
 The practical advantage of our approximations (they require approximating either a one-dimensional integral or an eigenvalue of an integral operator) is their relative simplicity as
to compute  the Glaz approximations one needs to numerically approximate $L+1$
and  $2L+1$   dimensional integrals. This is not an easy task even taking into account the fact of  existence of a sophisticated software; see references in Section~\ref{Glaz_ref}.

The paper is structured as follows. In Section~\ref{Reform} we reformulate the problem, state the Glaz approximation and discuss how to approximate our discrete-time problem with a continuous-time problem. In Section~\ref{Diffusion_section} we provide exact formulas for the first-passage probabilities (in the continuous-time setup) due to L. Shepp \cite{Shepp71} and give their alternative representation which will be crucial for deriving some of our approximations.   In Section~\ref{Corrected_shepp_section} we adapt the methodology of D. Siegmund to correct Shepp's formulas for discrete time and define  a version of the Glaz approximation which we will call Glaz-Shepp-Siegmund approximation.
In Section~\ref{GL_approx} we develop continuous-time approximations based on approximating eigenvalues of integral operators and subsequently correct them for discrete time. In Sections~\ref{CSA_sim_study} and \ref{sim_study2} we present results of large-scale simulation studies evaluating the performance of the considered approximations (also, in the cases when the original r.v. $\varepsilon_j$ are not normal and the weights in the moving window are not equal).
In Section~\ref{ARL_section}, we develop an approximation for ARL and compare its accuracy to the one developed in \cite{Glaz2012}.

\section{Boundary crossing probabilities: discrete and continuous time}\label{Reform}
\subsection{Standardisation of the moving sums}

The first two moments of $S_{n,L}$ are
$
\E S_{n,L}= \theta L$
and
$
{\rm var}(S_{n,L})=\displaystyle \sigma^2 L.
$
Define
\be
\label{eq:def-xi}
\xi_{n,L}:= \frac {S_{n,L}-  \E  S_{n,L}}
{\sqrt{{\rm var}(S_{n,L})}}=
\frac { S_{n,L}- \theta L}
{ \sigma\sqrt{  L    }}  \, ,\;\;\mbox{$n=0,1,\ldots\, ,$}
\ee
which are the standardized versions of $S_{n,L}$.
All r.v. $\xi_{n,L}$ are $N(0,1)$; that is, they have the probability density function and c.d.f.
\be
\label{eq:phi}
\varphi(x):=\frac{1}{\sqrt{2\pi}}e^{-x^2/2}\,,\; \; \Phi(t):= \int_{-\infty}^t \varphi(x) dx\, .
\ee
Unlike the original r.v. $\varepsilon_i$, the r.v.  $\xi_{0,L}, \xi_{1,L}, \ldots$  are correlated so that for all $k=0,1,\ldots$ we have
 ${\rm Corr}(\xi_{0,L} , \xi_{k,L})={\rm Corr}(\xi_{n,L},\xi_{n+k,L})$ and
\be
\label{eq:correlation}
{\rm Corr}(\xi_{n,L},\xi_{n+k,L})= \max\{0,1-{k}/{L}\}=
\left\{
  \begin{array}{cl}
    1-{k}/{L} & \;\;\;\;{\rm for}\; 0\le k \le L \\
    0 & \;\;\;\;{\rm for}\; k > L\, .
  \end{array}
\right.
\ee
Proof of \eqref{eq:correlation} is straightforward, {see \cite[Lemma 1]{AandZ2019}.}

Set  $T=M/L$ and
\be
\label{H_h}
h= \frac{H - \theta L}{\sigma \sqrt{L}}
\;\mbox{ so that }\;H= \theta L +
 \sigma h\sqrt{L} \,  .
\ee

Define  the BCP for the sequence of r.v. $\xi_{0,L}, \xi_{1,L}, \ldots$:
\be
 { P}_{L}(T,h) := \text{Pr}\left(\max_{n=0,1,\ldots,TL} \xi_{n,L} \geq  h \right)\, .
 \label{eq:prob-xi}
\ee
From \eqref{eq:def-xi} and \eqref{H_h}, the BCPs ${\cal P}_{\s}(M,H,L) $ and ${ P}_{L}(T,h)$ are equal:
\bea
{\cal P}_{\s}(M,H,L) = { P}_{L}(T,h)\,  \mbox{ for any }  H,L  \mbox{ and } T=M/L\,.
\eea
Note also that $ { P}_{L}(T,h) = 1-F_L(T,h)$, where
\be
 F_L(T,h) = \text{Pr}\left(\max_{n=0,1,\ldots,TL} \xi_{n,L} < h \right)\, .
 \label{eq:prob-xiF}
\ee
In accordance with the terminology of \cite{Shepp71} and \cite{slepian1961first} we shall call $ F_L(T,h)$ `first-passage probability'.
 In the following sections, we derive approximations for \eqref{eq:prob-xi}. These approximations will be
  based on approximating the sequence   of r.v. $\{\xi_{0,L}, \xi_{1,L}, \ldots,\xi_{M,L}\}$ by a continuous-time
random process and subsequently correcting the obtained approximations  for discreteness. Before doing this, we formulate the approximation which is currently the state-of-the-art.

\subsection{Glaz approximation for ${ P}_{L}(T,h)$}\label{Glaz_ref}
The approximation for the BCP ${ P}_{L}(T,h)$ developed in \cite{Glaz_old,Glaz2012, wang2014variable,wang2014multiple} and discussed in the introduction is as follows.\\

\noindent{\bf Approximation 1.} {\it (Glaz approximation) For $T \ge 2$,
\be
\label{eq:Glaz}
{ P}_{L}(T,h) \simeq 1- F_L(2,h)  \left[
 \frac{F_L(2,h) }{F_L(1,h) }
  \right]^{T-2} \, ,
\ee
where to approximate the first-passage probabilities  $F_L(1,h) $ and $F_L(2,h)$, which are $L+1$
and  $2L+1$   dimensional integrals respectively, it is advised to use the   so-called `GenzBretz' algorithm for numerical evaluation of multivariate normal probabilities; see  \cite{genz2009computation,GenzR}.
}\\

  Unless $h$ is large (say, $h>3$),   Approximation 1 is  very accurate. However,  its computational cost is also high, especially for large $L$.
   Moreover, the main option in the `GenzBretz' package requires the use of Monte-Carlo simulations so that
  for reliable estimation of high-dimensional integrals one needs to make a lot of averaging; see Section~\ref{sim_study_sub} and \ref{ARL_section}  for more discussion on these issues.

\subsection{{Continuous-time (diffusion) approximation}}
\label{sec:cont_time}

For the purpose of approximating the BCP ${ P}_{L}(T,h) $,    we replace the discrete-time process $\xi_{0,L}, \ldots,$ $ \xi_{M,L}$ with a continuous process $S(t)$, $t\in [0,T]$, where  $T=M/L$ (we will then correct the corresponding first-passage probabilities for discreteness).
    We do this as follows.

Set $\Delta= 1/L$ and define
$
t_n=n \Delta \in [0, { T}]
\; n=0,1,\ldots ,{M}.
$
Define a piece-wise linear  continuous-time process
${S_{L}(t)},$ $t \in [0,T]:$
\bea
{S_{L}(t)}=\!\frac1{\Delta} \left[
(t_{n}-t )\xi_{n-1,L}\! +\! (t-t_{n-1}) \xi_{n,L} \right] \;\;\;{\rm for}
 \;\;t \in [t_{n-1},t_n],\; n=1,\dots,{M}.\;
\eea
By construction, the process ${S_{L}(t)}$ is such that
${{S_{L}(t_n)}}=\xi_{n,L} \; {\rm for } \; n=0,\ldots,{M}$.
Also we have that ${S_{L}(t)}$ is  a second-order stationary process in the sense that
$ \E {S_{L}(t)},\,$ ${\rm var}({S_{L}(t)})$
and the autocorrelation  function
$R^{(L)}(t,t+k\Delta)={\rm Corr}({S_{L}(t)},{S_{L}(t+k\Delta)})$
do not depend on $t$.

\begin{lemma}
\label{Durbin} Assume $L \to \infty$.
The limiting process $S(t)$ = $\lim_{L \rightarrow \infty}{S_{L}(t)}$, where $t \in [0,T]$,
is a Gaussian second-order
stationary process with marginal distribution $S(t) \sim N(0,1)$ for all $t \in [0,T]$ and autocorrelation function $R(t,t+s)=R(s) = \max\{0,\;1\!-\!|s|\}\, $.
\end{lemma}
This lemma is a simple consequence of \eqref{eq:correlation}.

\subsection{Diffusion approximations: definition and their role in this study}

The above approximation of a discrete-time process $\{\xi_{0,L}, \xi_{1,L}, \ldots, \xi_{M,L}\}$ with a continuous process $S(t),\, t\in [0,T]$, allows us to approximate the BCP ${ P}_{L}(T,h)$  by a continuous-time analogue as follows.

By the definition of a diffusion approximation, the BCP
$
{ P}_{L}(T,h)
$
is approximated by
\be
{P}_{}(T,h) \!
:=\! {\rm Pr}\left \{\max_{0\leq t\leq { T}} S(t) \geq
h\right\}\, .
\label{eq:first_pass_prob}
\ee

Note that approximating the discrete process of moving sums by 
 a continuous-time process $S(t)$ and subsequent  approximation of the BCP ${ P}_{L}(T,h)$  by ${P}_{}(T,h)$ 
 is by no means new. This has been done, in particular, in \cite{Haiman}.

We will call  \eqref{eq:first_pass_prob} and any approximation to  \eqref{eq:first_pass_prob}, which does not involve the knowledge of $L$, `diffusion approximation'. These
 approximations can be greatly improved with the help of the methodology developed  by D.Siegmund and adapted to our setup in Section~\ref{Corrected_shepp_section}. The importance of the discrete-time correction is illustrated by Figures~\ref{L_5_CSA} and \ref{Importance of L}, where for a fixed $h$ and $T$ we can see a  significant difference in values of the BCPs
$
{ P}_{L}(T,h)
$
for different values of $L$. As seen from Figure~\ref{Importance of L}, even for very large $L=1000$, the discrete-time correction is still needed.
Hence we are not recommending to use any approximation for $
{ P}(T,h)
$
(including rather sophisticated ones like the one developed in \cite{Haiman})
as an approximation for $
{ P}_{L}(T,h)
$. In the next section we will discuss a diffusion approximation that, after correcting for discrete time, will be a cornerstone for all approximations developed in this paper.

In what follows, it will also be convenient to use the first-passage probability
 \bea
 {F}_{}(T,h) = {\rm Pr}\left \{\max_{0\leq t\leq { T}} S(t) <
h\right\} =  1- {P}_{}(T,h)\, .
\eea
Since $\xi_{0,L}=S(0) \sim N(0,1)$, we have ${F}_{}(0,h)=1-{P}_{}(0,h) = \Phi(h)   $.

\section{Exact formulas for the first-passage probabilities in the continuous-time case }\label{Diffusion_section}

\subsection{Shepp's formulas}
\label{sec:Shepp}
Define the conditional first-passage probability
\be \label{F_conditional}
{ F}_{}(T,h\, |\,  x) := {\rm Pr}\Big\{ S(t) \! < \!
h\;{\rm for\; all\; } t\!\in \! [0,\!{ T}]  \, | \, S(0) = x \Big\} \, . \;\;
\ee
Since ${ F}_{}(T,h\, |\,  x) = 0$ for $x>h$, for the unconditional first-passage probability $ {F}_{}(T,h)$ we have
 $ {F}_{}(T,h) = \int_{-\infty}^{h}{ F}_{}(T,h\, |\,  x)\varphi(x)dx$.

The result of  \cite[p.949]{Shepp71} states than if $T=n$ is a positive integer then
\begin{equation}\label{shepp_form}
{ F}_{}(n,h\, |\,  x) = \frac{1}{\varphi(x)} \int_{D_x} \det[\varphi(y_i - y_{j+1} + h)]^n_{i,j=0} \, dy_2\ldots dy_{n+1}\,
\end{equation} where $y_0= 0, y_1=h-x,$
$
D_x =  \{y_2, \dots , y_{n+1} \given h-x < y_2 < y_3 < \ldots < y_{n+1}   \}
$.
For non-integer $T\ge1$,  the exact formula for $
{ F}_{ }(T,h\, |\,  x)
$ is even more complex (the integral has the dimension $\lceil 2T \rceil $) and completely impractical for computing ${P}_{}(T,h)$ with $T>2$,  see \cite[p.950]{Shepp71}.

For $n=1$, we obtain
\be \label{F1_form}
{F}_{}(1,h) &=&  \int_{-\infty}^{h}\int_{-x-h}^{\infty} \det\begin{bmatrix}
    \varphi(x)      &  \varphi(-x_2\!-\!h) \\
    \varphi(h)      & \varphi(-x\!-\!x_2)\\
\end{bmatrix} dx_2dx\,
=   \Phi(h)^2 - \varphi({h})[h\Phi(h)+\varphi(h)].
\ee
For $n=2$, \eqref{shepp_form} yields
\be
{F}_{}(2,h)\!&=&\! \! \int_{-\infty}^{h}\int_{-x-a}^{\infty}\int_{x_2-a}^{\infty}\! \det\begin{bmatrix}
    \varphi(x)      &  \varphi(-x_2\!-\!a) &  \varphi(-x_3\!-\!2a) \\
    \varphi(a)      & \varphi(-x\!-\!x_2)&   \varphi(-x\!-\!a\!-\!x_3) \\
    \varphi(x_2\!+\!2a\!+\!x)      & \varphi(a)&  \varphi(x_2\!-\!x_3) \\
\end{bmatrix}\! \! dx_3dx_2dx . \label{F2_form}
\ee
The three-dimensional integral in \eqref{F2_form} can be reduced to a one-dimensional, see \eqref{F_Corrected_shepp_explicit_2} below with $h_L=h$.



\subsection{An alternative  representation of the Shepp's  formula \eqref{shepp_form}}
\label{sec:altern}

 Set $s_i=h+y_i-y_{i+1}$ ($i=0,1, \ldots, n$)   with $s_0=x$, $y_0= 0, y_1=h-x$.
It follows from Shepp's proof of \eqref{shepp_form} that  $s_0, s_1, \ldots, s_n$ have the meaning of the values of the process $S(t)$ at the times $t=0,1, \ldots, n$:
$S(i)=s_i$ ($i=0,1, \ldots, n$). The range of the variables $s_i$ is $(-\infty,h)$.
Changing the variables in \eqref{shepp_form}, we obtain
\begin{equation}\label{shepp_form2}
{ F}_{}(n,h\, |\,  x)  = \frac{1}{\varphi(x)} \int_{-\infty}^h \ldots \int_{-\infty}^h  \det[\varphi(s_i + a_{i,j})]^{n}_{i,j=0} \, ds_1\ldots ds_{n}\,,
\end{equation}
where
\bea
a_{i,j}=y_{i+1} \!- \!y_{j+1}=\left\{ \begin{array}{cl}
 0 & \text{   for } i=j\, \\
     (i-j)h \!- \!s_{j+1}\!-\!\ldots\!- \!s_{i+1} & \text{   for } i>j\,  \\
       (i-j)h+  s_{i+1}+\ldots+ s_{j}  & \text{   for } i<j\, .

       \end{array} \right.
\eea

\subsection{Joint density for the values $\{S(i)\}$ and  associated transition densities}
From \eqref{shepp_form2}, we obtain the following expression for
the joint probability density function for the values $S(0),S(1), \ldots, S(n)$ under the condition  $S(t)<h$ for all $t\in [0,n]$:
\begin{equation}\label{shepp_form7}
p(s_0,s_1,\ldots s_{n})= \frac{1}{\varphi(s_0) { F}_{}(n,h\, |\,  s_0) }  \det[\varphi(s_i +a_{i,j})]^{n}_{i,j=0} \, .
\end{equation}
From this  formula, we can derive   the transition density from $s_0\!=\!x$ to $s_n$ conditionally
$S(t)\!<\!h,$ $\forall t\in [0,n]$:
 \begin{equation}\label{shepp_form9}
q^{(0,n)}_h(x\to s_n ) = \frac{1}{\varphi(x)} \int_{-\infty}^h \ldots \int_{-\infty}^h  \det[\varphi(s_i + a_{i,j})]^{n}_{i,j=0} \, ds_1\ldots ds_{n-1}\, .
\end{equation}
For this transition density, $ \int_{-\infty}^h   q^{(0,n)}_h(x\to z ) dz = { F}_{}(n,h\, |\,  x) $. Moreover, since $S(0)\sim N(0,1)$, the non-normalized density of $S(n)$ under the condition $S(t)<h$ for all $t\in[0,n]$ is
\be\label{p_n_density}
 {p}^{(0,n)}_h(z) :=  \int_{-\infty}^{h}q^{(0,n )}_h(x\to z )\varphi(x) dx
\ee
with $z<h$ and $ \int_{-\infty}^h    {p}^{(0,n)}(z) dz = { F}_{}(n,h) $.
In the case  $n=1$,  \eqref{shepp_form9} gives
 \begin{equation}\label{shepp_form3}
q^{(0,1)}_h(x\to z ) = \frac{1}{\varphi(x)}   \det  \left(
                                                      \begin{array}{cc}
                                                       \varphi(x)  & \varphi(x\!-\!h\!+\!z) \\
                                                        \varphi(h) & \varphi(z) \\
                                                      \end{array}
                                                    \right)= \varphi(z)\left[1 - e^{-(h-z)(h-x)} \right]\,, \,\,  z=s_1<h.
\end{equation}
From this and \eqref{p_n_density} we get
 \bea
 {p}^{(0,1)}_h(z) = \int_{-\infty}^{h}q^{(0,1 )}_h(x\to z )\varphi(x) dx = \Phi(h)\varphi(z) - \Phi(z)\varphi(h)
 \eea
with $z<h$ and $ \int_{-\infty}^h    {p}^{(0,n)}(z) dz = { F}_{}(1,h) $.

Rather than just recovering the transition density from $s_0=x$ to $s_n$, we can also use \eqref{shepp_form7} and \eqref{p_n_density} to  obtain the transition density from $x=s_j$ to $z=s_n$, $0< j<n$, under the condition $S(t)<h$ for all $t\in[0,n]$:
\be
\label{0shepp_form9}
{q}_h^{(j,n)}(x\!\to\! z )\!=\!\frac{1}{ {p}^{(0,j)}_h(z)}\int_{-\infty}^h \!\ldots\! \int_{-\infty}^h \!\det[\varphi(s_i\! +\! a_{i,j})]^{n}_{i,j=0} \, ds_0ds_1\ldots ds_{j-1}ds_{j+1}\ldots ds_{n-1},\;\;\;\;\;
\ee
where $s_j=x$ and $s_n=z$.
For $j=1$ and $n=2$
we obtain  the transition density from $x=s_1$ to $z=s_2$ under the condition  $S(t)<h$ for all $t\in [0,2]$:
\be\label{q_x_z}
{q}_h^{(1,2)}(x\to z )& =&
   \frac{1}{ {p}^{(0,1)}_h(z)}\int_{-\infty}^h
   \det  \left(
                                                      \begin{array}{ccc}
                                                       \varphi(s_0)  & \varphi(s_0\!-\!h\!+\!x) &\varphi(s_0\!-\!2h\!+\!x\!+\!z) \\
                                                        \varphi(h) & \varphi(x)& \varphi(\!x\!+\!z\!-\!h) \\
                                                        \varphi(2h\!-\!x) & \varphi(h)& \varphi(z)
                                                      \end{array}
                                                    \right) ds_0
\nonumber \\ &=&
\frac{1}{\Phi(h)\varphi(x) - \Phi(x)\varphi(h)} \det  \left(
                                                      \begin{array}{ccc}
                                                       \Phi(h)  & \Phi(x) &\Phi(\!x\!+\!z-h) \\
                                                        \varphi(h) & \varphi(x)& \varphi(\!x\!+\!z\!-\!h) \\
                                                        \varphi(2h\!-\!x) & \varphi(h)& \varphi(z)
                                                      \end{array}
                                                    \right) \, .
\ee

\section{Correcting Shepp's formula \eqref{shepp_form} for discrete time}
\label{Corrected_shepp_section}

\subsection{Rewriting  \eqref{shepp_form} in terms of the Brownian motion }

Let $W(t)$ be the standard Brownian Motion process on $[0, \infty)$ with $W(0)=0$ and
$
	\E W(t) W(s) =\min(t, s).
$ Recall the conditional probability ${ F}_{}(T,h\, |\,  x)$ defined in \eqref{F_conditional}. Suppose $T\ge1$ is an integer and define the event
\bea
{ \rm \Omega} &=& \{W(t)  < W(t+1)+h< W(t+2)+2h< \cdots < W(t+T)+Th , \,\, \forall \,\,\,0 \le t \le 1 \}\\
 &=& \{W(t) -W(t+1)<h,\ldots, W(t+T-1) -W(t+ T)<h , \,\, \forall \,\,\,0 \le t \le 1 \}.
\eea
If $W(i)=x_i$, $i=0,1,\ldots, T+1$, we obtain from \cite[p.948]{Shepp71}
\begin{eqnarray}\label{condition}
{ F}_{}(T,h\, |\,  x) = \int \cdots \int {\rm Pr}\{ { \rm \, \Omega} \, \big| \, W(i) =x_i,\,\, i=0,1,2,\ldots,T+1,\,\,\,  W(0)=0,\,\, W(0)-W(1) = x\} \nonumber \\ \times \, {\rm Pr} \{ W(i) \in dx_i, \,\, i=0,1,2,\ldots,T+1,\,\, \big|\,\,   W(0)=0,\,\, W(0)-W(1) = x\}. \nonumber \\
\end{eqnarray}
It follows from the proof of \eqref{shepp_form} that to correct \eqref{condition} for discrete time, one must correct the following probability for discrete time
\begin{eqnarray}
&&{\rm Pr}\{ { \rm \Omega}  \,\, \big| \,\, W(i) =x_i,\,\, i=0,1,2,\ldots,T+1,\,\,\,  W(0)=0,\,\, W(0)-W(1) = x\} \nonumber\\
&=&{\rm Pr}\{ \sqrt{2}\, W_1(t) <h , \ldots,  \sqrt{2}\, W_T(t) <h ,  \,\, \forall \,\,\,0 \le t \le 1 \,\, \big| \,\, W(i) =x_i,\,\, i=0,1,2,\ldots,T+1,\,\,\, \nonumber \\
  &&\qquad\qquad\qquad\qquad\qquad\qquad\qquad\qquad\qquad\qquad\qquad\quad  W(0)=0,\,\, W(0)-W(1) = x\}\label{BM_condition}
\end{eqnarray}
where $W_i(t) =\frac{\sqrt{2}}{2} [W(t+i-1) -W(t+i)] $, $i=1,2,\ldots,T$. Due to the conditioning on the rhs of \eqref{BM_condition},  the processes $W_i(t)$ can be treated as independent Brownian motion processes. Therefore, the independent increments of the Brownian motion means correcting formula \eqref{shepp_form} for discrete time is equivalent  to correcting the probability ${\rm Pr}(\sqrt{2}\, W(t) <h, \,\, \forall \,\,\,0 \le t \le 1 \,\, )$ for discrete time.

\subsection{Discrete-time correction for the BCP of cumulative sums.}\label{Siegmund_discrete_time}

Let $X_1,X_2, \ldots $ be i.i.d. $N(0,1)$ r.v's and set $Y_n = X_1+X_2+ \ldots + X_n$.
Consider the sequence of cumulative sums $\{Y_n\}$ and define the stopping time
$
\tau_{Y,a,b} = \inf \{ n \ge 1: Y_n \ge a + bn \}
$
for $a > 0$ and $b \in \mathbb{R}$. Consider the problem of evaluating
\begin{equation}
\label{51}
{\rm Pr}(\tau_{Y,a,b}\le N) = {\rm Pr}(Y_n \ge a+bn \text{ for at least one } n \in \{1,2,\ldots N \}).
\end{equation}

Exact evaluation of \eqref{51} is difficult even if $N$ is not very large but it was accurately approximated by D.Siegmund see e.g. \cite[p.19]{Sieg_paper}. Let $W(t)$ be the standard Brownian Motion process on $[0, \infty)$. For $a>0$ and $b \in \mathbb{R}$, define
$
\tau_{W,a,b} = \inf \{ t: W(t) \ge a + bt \}
$
so that
\begin{equation}\label{D_approx}
{\rm Pr}(\tau_{W,a,b}\le N) = P_W(N, a+bt) := {\rm Pr}\left\{W(t) > a+bt \text{ for at least one }  t\in[0,N] \right\}\,.
\end{equation}

In \cite{Sieg_paper}, \eqref{D_approx} was used to approximate \eqref{51} after translating the barrier $a+bt$ by a suitable scalar $\rho \ge 0$. Specifically, the following approximation has been constructed: \bea
P(\tau_{Y,a,b}\le N) \cong P_W(N, (a+\rho)+bt)\, ,
\eea
where the constant $\rho$ approximates the expected excess of the process $\{Y_n\}$ over the barrier  $a+bt$.
From \cite[p. 225]{Sieg_book}
\begin{equation}
\label{D_rho}
\rho =  - \pi^{-1}\int_{0}^{\infty}\lambda^{-2}  \log\{2(1-\exp(-\lambda^2/2))/ \lambda^2 \} \, d\lambda\, \simeq 0.582597.
\end{equation}

\subsection{Discretised Brownian motion }\label{Q_h_correction}

 Define $\epsilon$ = $1/L$ and let
$
t^\prime_n= n\epsilon \in [0, { 1}],
$
$ n=0,1,\ldots ,{L}.$
Let $X_1,X_2, \ldots $ be i.i.d. $N(0,1)$ r.v's and set
$
W(t^\prime_n) = \sqrt{\epsilon}\sum_{i=1}^{{n}}X_i.
$
For $a>0$ define the stopping time
\begin{equation}\label{Brown_stopping}
\tau_{W,a,b} = \inf \{t^\prime_n: \sqrt{2}W(t^\prime_n) \ge a  \}
\end{equation}
and consider the problem of approximating
\begin{equation}
\label{Problem} \!\!\!\!\!
\!{\rm Pr}(\tau_{W,a,b} > 1) = {\rm Pr}\bigg(\sqrt{2}W(t^\prime_n)  < a\text{ for all } t^\prime_n \in \{0, \epsilon, \ldots,L\epsilon=1 \}\bigg).
\end{equation}
As $L \rightarrow \infty$, the piecewise linear continuous-time process $W^\epsilon(t)$, $t \in [0,1]$, defined by:
\begin{equation*}
W^\epsilon (t)\!:=\!\frac1{\epsilon} \left[
(t^\prime_{n}-t )W(t^\prime_{n-1})\! +\! (t-t^\prime_{n-1}) W(t^\prime_{n})\right] \;\;\;{\rm for}
 \;\;t \in [t^\prime_{n-1},t^\prime_n],\; n=1,\dots,{L},\;
\end{equation*}
converges to $W(t)$ on $[0,1]$ as so we can refer to $W(t^\prime_n)$ as discretised Brownian motion.  We make the following connection between $\sqrt{2}W(t^{\prime}_n)$ and the random walk $Y_n$:
\begin{equation*}
\sqrt{2}W(t^\prime_n)=\sqrt{2\epsilon} \,Y_n = \frac{\sqrt{2}}{\sqrt{L}}{Y_n}\,\,, n=1,2,\ldots M.
\end{equation*}
Then by using \eqref{D_rho}, we approximate the expected excess over the boundary for the process $\sqrt{2}W(t^\prime_n)$  by
\bea
\omega_L := \frac{{0.82}}{\sqrt{L}} \simeq   \frac{\sqrt{2} \rho}{\sqrt{L}} \, .
\eea
We have deliberately rounded the value $\sqrt{2} \rho \simeq 0.8239...$ to $0.82$ as for  small $h$ and small $L$ it provides marginally better approximation  \eqref{shepp_correction}.

\subsection{Corrected version of \eqref{shepp_form}} \label{sec:corr}
Set ${h}_L = h+\omega_L$. To correct \eqref{shepp_form} for discrete time we substitute the barrier $h$ with $h_L$. From this and the relation $ {F}_{}(T,h) = \int_{-\infty}^{h}{ F}_{}(T,h\, |\,  x)\varphi(x)dx$, the discrete-time corrected form of ${ F}_{}(T,h) $  is
\be\label{correct_F_T}
 { F}_{}(T,{h},h_L) &:=& \int_{\infty}^{h}{ F}_{}(T,{h_L}\, |\,  x) \varphi(x)dx
  =\int_{-\infty}^{h}\int_{D_x} \det[\varphi(y_i - y_{j+1} + {h_L})]^{T}_{i,j=0} \, dy_2\ldots dy_{T+1}\,dx,\nonumber\\
  &&
\ee
where $ y_0= 0, y_1={h_L}-x,$ and
$
D_x =  \{y_2, \dots , y_{T+1} \given {h_L}-x < y_2 < y_3 < \ldots < y_{T+1}   \}.
$

\subsection{A generic approximation involving corrected Shepp's formula}\label{CSA_sec}

\noindent{\bf Approximation 2.} \textit{For integral $T\ge 1$,  the discrete-time correction for the BCP \eqref{eq:prob-xi} is
\be\label{shepp_correction}
{ P}_{L}(T,h) \cong {P}_{}(T,{h},h_L):=
1-  { F}_{L}(T,{h},h_L),
\ee
where ${ F}_{L}(T,{h},h_L)$ is given in \eqref{correct_F_T}. }\\

Whilst Approximation 2 is very accurate (see the next subsection), computation of ${P}_{}(T,{h},h_L)$ requires numerical evaluation of a $T+1$ dimensional integral which is impractical for large $T$. To overcome this, in Section~\ref{corrected_diffusions} we develop approximations that can be easily used for any $T>0$ (which is not necessarily integer).

\subsection{Particular cases: $T=1$ and $T=2$}\label{sim_study1}\label{T1T2}

For $T=1$, evaluation of \eqref{correct_F_T} yields
\begin{equation}\label{F_Corrected_shepp_explicit}
{F}_{}(1,{h},h_L) = \Phi(h)\Phi({h_L}) - \varphi({h_L})[h\Phi(h)+\varphi(h)]\, .
\end{equation}
In our previous work  \cite{AandZ2019} we have derived approximations
$\hat{ P}_{L}(T,h)$ for the BCP ${ P}_{L}(T,h)$ with $T\leq 1$.
The approximations $\hat{ P}_{L}(T,h)$ developed in \cite{AandZ2019} are also discrete-time corrections of the
continuous-time probabilities ${ P}(T,h)$ but they are based almost exclusively on   the fact that the process $S(t)$ is conditionally Markov on the interval $t\in [0,1]$;  hence the technique of \cite{AandZ2019} cannot be extended for intervals $t \in [0,T]$ with $T>1$.
The approximation $\hat{ P}_{L}(1,h)$ of \cite{AandZ2019} is
different from ${P}_{}(1,{h},h_L) =1-{F}_{}(1,{h},h_L) $ of \eqref{F_Corrected_shepp_explicit}.
It appears that   $\hat{ P}_{L}(1,h)$ is more complicated and less accurate approximation than ${P}_{}(1,{h},h_L)$.

For $T=2$,  \eqref{correct_F_T} can be expressed (after some manipulations) as follows:
\be \label{F_Corrected_shepp_explicit_2}
\!\!\!\!\!{F}_{}(2,{h},h_L)\! &=&\!
  \frac{\varphi^2(h_L)}2 \!  \left[ ({h}^{2}\!-\!1\!+\!\sqrt{\pi}h)\Phi \left( h \right) \!+\!(h\!+\!\sqrt{\pi})\varphi \left( h \right)\!
  \right] \!\! -\!\varphi \left( h_L \right) \Phi \left(  h_L \right)  \left[  \left( h\!+\!
h_L  \right) \Phi \left( h  \right)\! +\!\varphi \left( h \right)  \right]\nonumber \\
 &+&\!\Phi \left( h \right)   \Phi^2({ h_L} ) +\!\! \int_{0}^{\infty} \!\Phi(h\!-\!y)\left[\varphi(h_L+y)\Phi(h_L-y)\!-\!\sqrt{\pi}\varphi^2(h_L)\Phi(\sqrt{2}y) \, \right] \!dy.\,
  \ee
Only a one-dimensional integral has to be numerically evaluated for computing  ${F}_{}(2,{h},h_L)$.

\subsection{Simulation study}\label{sim_study1}\label{CSA_sim_study}

\begin{figure}[h]
\begin{center}
 \includegraphics[width=0.5\textwidth]{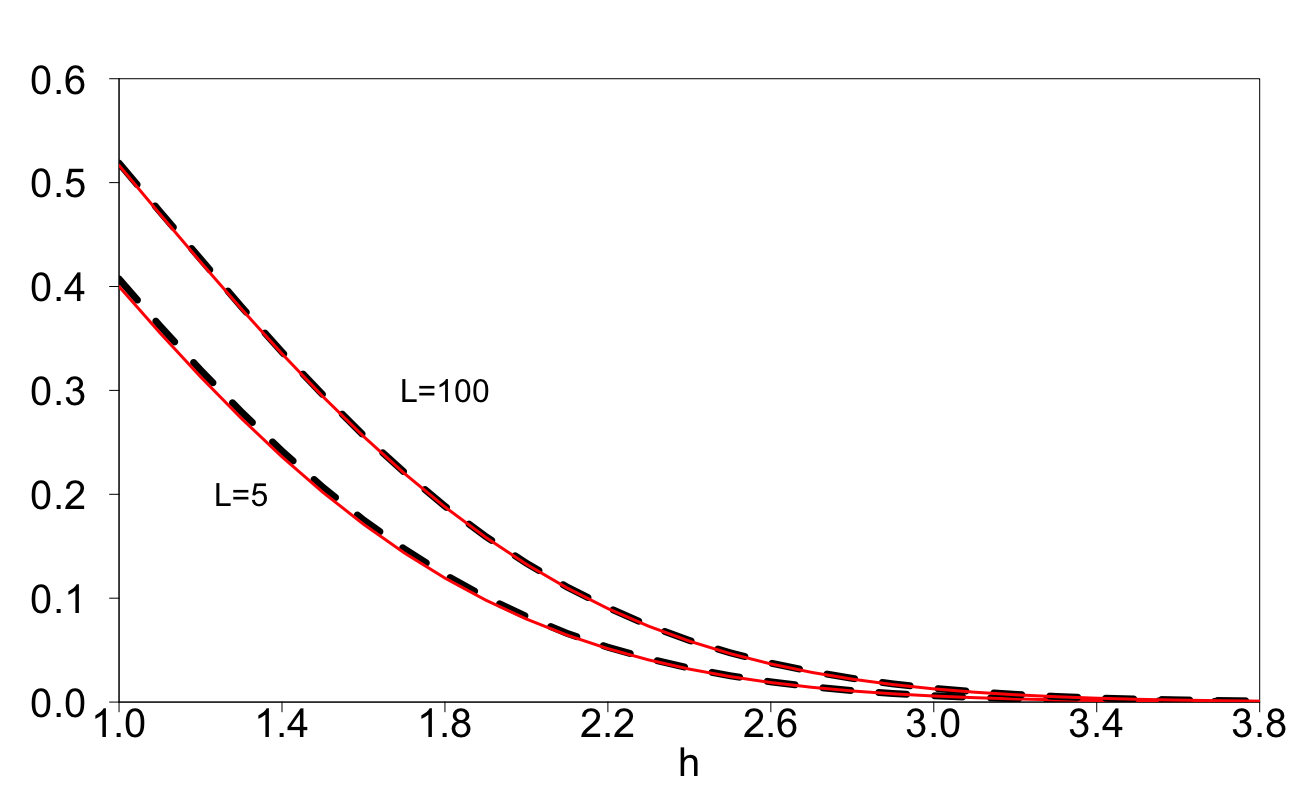}$\;\;$\includegraphics[width=0.5\textwidth]{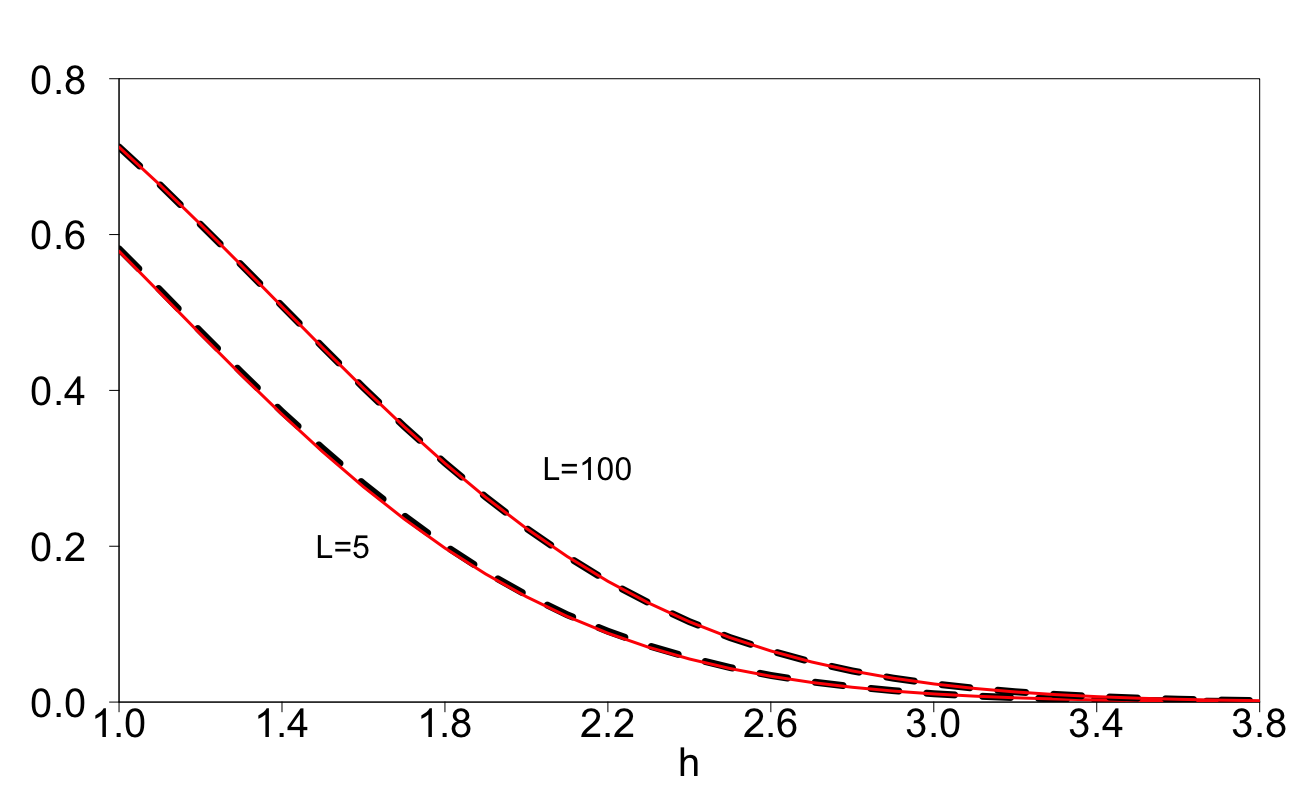}
\end{center}
\caption{Empirical probabilities of reaching the barrier $h$ (dashed black) and corresponding versions of Approximation 2 (solid red).
Left: $T=1$ with (a) $L=M=5$ and (b) $L=M=100$. Right: $T=2$ with (a) $L=5$, $M=10$ and (b)  $L=100$, $M=200$ . }
\label{L_5_CSA}
\end{figure}

\begin{figure}[h]
\begin{center}
\includegraphics[width=0.5\textwidth]{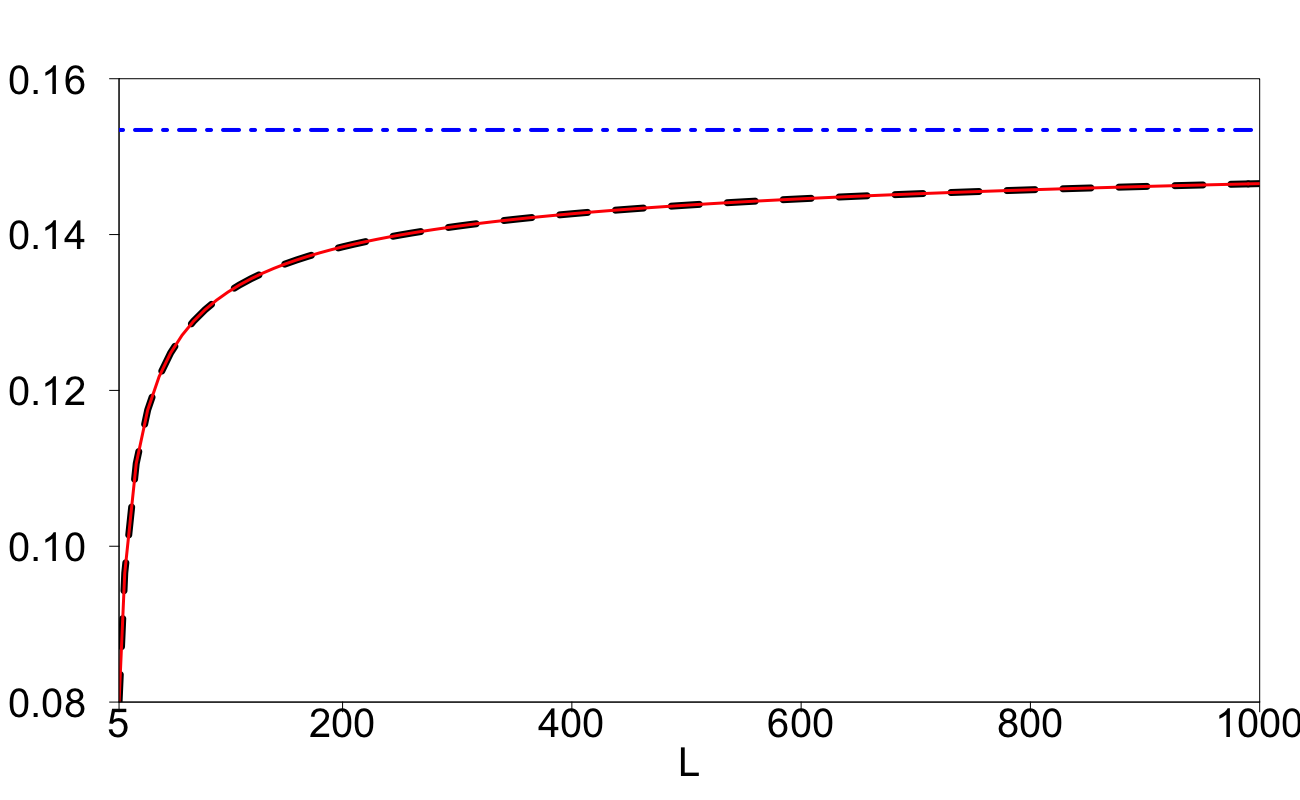}$\;\;$\includegraphics[width=0.5\textwidth]{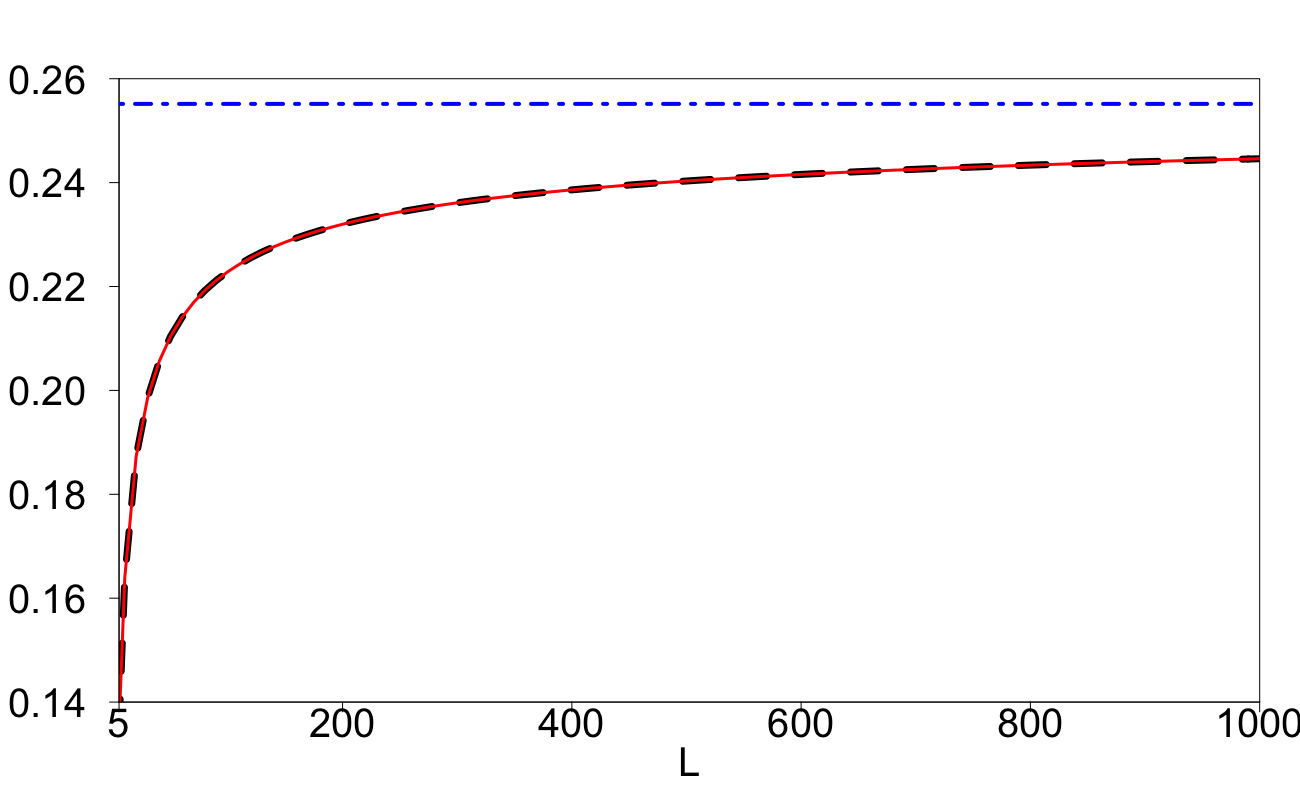}
\end{center}
\caption{Empirical probabilities of reaching the barrier $h=2$ as a function of $L$ (dashed black), uncorrected diffusion approximation ${P}_{}(T,{2})$ (dot-dashed blue) and corresponding version of ${ P}_{L}(T,h)$, which is  Approximation 2 (solid red).
Left: $M\!=\!L$ ($T\!=\!1$). Right: $M\!=\!2L$ ($T\!=\!2$). }
\label{Importance of L}
\end{figure}

In this section, we assess the quality of the approximations \eqref{F_Corrected_shepp_explicit} and \eqref{F_Corrected_shepp_explicit_2} as well as  the sensitivity of the BCP ${ P}_{L}(T,h)$ to the value of $L$.  In Figures~\ref{L_5_CSA} and \ref{Importance of L}, the black dashed line corresponds to the empirical values of the BCP ${ P}_{L}(T,h)$ (for $T=M/L=1,2$) computed from 100\,000 simulations  with different values of $L$ and $M$ (for given $L$ and $M$, we simulate $L+M$ normal random variables 100\,000 times). The solid red line corresponds to Approximation~2. The axis are: the $x$-axis shows the value of the  barrier $h$ in Figure~\ref{L_5_CSA} and value of $L$ in Figure~\ref{Importance of L}; the $y$-axis denotes the probabilities of reaching the barrier. The graphs, therefore,  show the empirical probabilities of reaching the barrier $h$ (for the
dashed line) and values of considered approximations for these probabilities. From these graphs we can conclude that Approximation~2 is very accurate, at least for $T=1,2$. We can also conclude that the BCP ${ P}_{L}(T,h)$ is very sensitive to the value of $L$. From Figure~\ref{Importance of L} we can observe a counter-intuitive fact that even for very high value $L=1000$, the BCP ${P}_{L}(T,{h})$ is not even close to
$ {P}_{\infty}(T,{h}) = {P}_{}(T,{h}) $ from \eqref{eq:first_pass_prob}. This may be explained by the fact that for any fixed $T$ and $h$, the  inaccuracy  $|{ P}_{L}(T,h)- {P}_{}(T,{h})|$ decreases with the rate const$/\sqrt{L}$ as $L \to \infty$.


\subsection{The Glaz-Shepp-Siegmund approximation}

Combining \eqref{eq:Glaz} and  the  approximation \eqref{shepp_correction} for Shepp's formula \eqref{shepp_form}, we arrive at
the following approximation to which we suggest the name   `Glaz-Shepp-Siegmund approximation'.\\

\noindent{\bf Approximation 3.} \textit{For all $T >0, $ }
\be\label{GSSA_approx}
{ P}_{L}(T,h)\simeq 1- { F}_{}(2,{h},h_L) \cdot  \mu_L(h)^{T-2} \;\;\textit{with} \;
\mu_L(h)= \frac{{ F}_{}(2,{h},h_L) }{{ F}_{}(1,{h},h_L) }\, ,
\ee
\textit{where ${ F}_{}(1,{h},h_L) $ and  ${ F}_{}(2,{h},h_L)$ are defined in \eqref{F_Corrected_shepp_explicit} and  \eqref{F_Corrected_shepp_explicit_2} respectively.
}\\\\
Approximations 1 and 3 look  similar but  computing Approximation 1 is very hard and Approximation 3 is very easy (only a one-dimensional integral should be numerically computed).

\section{Approximations for the BCP ${ P}_{L}(T,h)$ through eigenvalues of integral operators }\label{GL_approx}

\subsection{Continuous time: approximations for $F(T,h)$}\label{one_step_sec}

Let $m$ be a positive integer, and  $q(x\to z )$ be the transition density $ {q}_h^{(m-1,m)}(x\to z )$ defined by \eqref{shepp_form3} for $m=1$ \eqref{q_x_z} for $m=2$ and
\eqref{0shepp_form9} for $m>2$.

Let us approximate the distributions of the values $s_i=S(i)$ for integral $i>m$ in the following way. Let   $p_{i}(x)$   be the density of $S(i) $ under the condition that  $S(t)$ does not reach  $h$ for $t \in [0,i]$. By ignoring the past values of $S(t)$ in $[0,i)$, the non-normalized density of
$S(i+1)$ under the conditions that $S(i)\sim p_{i}(x)$ and $S(t)$ does not reach  $h$ for $t \in [i,i+1]$ is
 \begin{equation}
\label{2.24}
\tilde{p}_{i+1}(x) = \int_{-\infty}^{h}{q}_h(x\to z )p_{i}(y)dy, \text{    for $x<h$}\,.
\end{equation}
 We can then define ${p}_{i+1}(x)=\tilde{p}_i(x)/c_{i}, \; x<h,$
where $c_{i}=\int_{-\infty}^{h}\tilde{p}_{i}(x)dx$.
We then replace formula (\ref{2.24}) with
\begin{equation}
\label{2.26}
\tilde{p}_i(x) = \int_{-\infty}^{h}{q}_h(x\to z )p(y)dy, \text{    for $x<h$},
\end{equation}
where $p(x)$ is an eigenfunction of the integral operator with kernel \eqref{shepp_form3} corresponding to the maximum eigenvalue $\lambda_m(h)$:
\be
\lambda_m(h) p(x) = \int_{-\infty}^{h}p(y){q}_h^{(m-1,m)}(x\to z )dy, \text{  } x<h\, . \label{2.27_1}
\ee
This eigenfunction $p(x)$ is a probability density on $(-\infty,h]$ with $p(x)> 0$ for all $x \in (-\infty,h)$ and
$
\int_{-\infty}^{h}p(x)dx = 1 \, .\label{2.28}
$
Moreover, the maximum eigenvalue $\lambda_m(h)$ of the operator with kernel $K(x,y)={q}_h^{(m-1,m)}(x\to z )$ is simple and positive.
The fact that such maximum eigenvalue $\lambda_m(h)$ is  simple and real (and hence positive) and the eigenfunction $p(x)$ can be chosen as  a probability density follows from the Ruelle-Krasnoselskii-Perron-Frobenius theory of bounded linear positive operators, see e.g. Theorem XIII.43 in \cite{ReedSimon}.

Using  (\ref{2.26}) and (\ref{2.27_1}), we  derive recursively:
$
{ F}(i+1,h) \simeq   F(i,h)\lambda_m(h)$ ($i=m,m+1,\ldots$).
By induction,  for any integer $T\geq m$ we then have
\be
\label{appr_main_1}
 {F}_{}(T,h) \simeq \,
                                           {F}_{}(m,h)\cdot \left[ \lambda_m(h)\right]^{T-m} \,.
                                           \ee
The approximation \eqref{appr_main_1} can  be used for any $T>0$ which is not necessarily an integer. The most important particular cases of \eqref{appr_main_1} are with $m=1$ and $m=2$. In these two cases, the kernel ${q}_h^{(m-1,m)}(x\to z )$ and hence the approximation \eqref{appr_main_1}  will be corrected for discrete time in the next section.

\subsection{Correcting approximation \eqref{appr_main_1}  for discrete time}\label{corrected_diffusions}

To correct the approximation \eqref{appr_main_1}  for discrete time we need to correct: (a) the first-passage probability ${F}_{}(m,h)$ and (b)
 the kernel ${q}_h^{(m-1,m)}(x\to z )$. The discrete-time correction of ${F}_{}(m,h)$ can be done using  ${F}_{L}(m,h,h_L)$ from  \eqref{correct_F_T} so that what is left is to correct the kernel ${q}_h^{(m-1,m)}(x\to z )$ and hence $\lambda_m(h)$.

\subsubsection{Correcting the transition kernels for discrete time}

As explained in Section \ref{Corrected_shepp_section}, to make a discrete-time correction in the Shepp's formula \eqref{shepp_form}  we need to replace the barrier $h$ with  ${h_L}=h+\omega_L$ in all places except for the upper bound for the initial value  $S(0)$. Therefore, using the notation of Section~\ref{sec:altern}, the joint probability density function for the values $S(0),S(1), \ldots, S(m)$ under the condition  $S(t)<h$ for all $t\in [0,m]$ corrected for discrete time is:
\begin{equation}\label{shepp_form_joint_density}
\hat{p}(s_0,s_1,\ldots s_{m})= \frac{1}{\varphi(s_0) { F}_{}(m,h\, |\,  s_0) }  \det[\varphi(s_i +\hat{a}_{i,j})]^{m}_{i,j=0} \,\
\end{equation}
with $-\infty<s_0<h$, $-\infty<s_j<h_L$ $(j=1,\ldots,m)$,
\bea
\hat{a}_{i,j}=y_{i+1} \!- \!y_{j+1}=\left\{ \begin{array}{cl}
 0 & \text{   for } i=j\, \\
     (i-j){h_L} \!- \!s_{j+1}\!-\!\ldots\!- \!s_{i+1} & \text{   for } i>j\,  \\
       (i-j){h_L}+  s_{i+1}+\ldots+ s_{j}  & \text{   for } i<j\,.
       \end{array} \right. \,
\eea

This gives us   the discrete-time corrected transition density from $s_0\!=\!x$ to $s_m$ conditionally
$S(t)\!<\!h,$ $\forall t\in [0,m]$:
 \begin{equation}\label{shepp_form9_correct}
q^{(0,m)}_{h_L}(x\to s_m ) = \frac{1}{\varphi(x)} \int_{-\infty}^{{h_L}} \ldots \int_{-\infty}^{{h_L}}  \det[\varphi(s_i + \hat{a}_{i,j})]^{m}_{i,j=0} \, ds_1\ldots ds_{m-1}\, ;
\end{equation}
which is exactly  \eqref{shepp_form9} with  $h_L$ is substituted for $h$.
In a particular case $m=1$, the corrected  transition density is
\be\label{Corrected_one_Step}
q^{(0,1)}_{h_L}(x\to s_1 )  =  \frac{1}{\varphi(x)}   \det  \left(
                                                      \begin{array}{cc}
                                                       \varphi(x)  & \varphi(x\!-\!{h_L}\!+\!s_1) \\
                                                        \varphi({h_L}) & \varphi(s_1) \\
                                                      \end{array}
                                                    \right)= \varphi(s_1)\left[1 - e^{-({h_L}-s_1)({h_L}-x)} \right]
\ee
 with $s_1<{h_L}$.

 Let us now make the discrete-time correction of the transition density $ {q}_h^{(1,2)}(x\to z )$.
 Denote by ${p}_{h,L}^{(0,1)}(z)$, $z<h,$  the non-normalized density of
 $S(1)$ under the condition $S(t)<h$ for all $t \in[0,1]$ corrected for discrete time; it satisfies $\int_{-\infty}^{h}{p}_{h,L}^{(0,1)}(z)dz ={ F}_{}(1,h,h_L) $. Using \eqref{Corrected_one_Step}, we obtain
 \bea
{p}_{h,L}^{(0,1)}(z) = \int_{-\infty}^{h}q^{(0,1)}_{h_L}(x\to z )\varphi(x) dx = \varphi(z)\Phi(h)-\varphi({h_L})\Phi(h-{h_L}+z) .
 \eea
From  \eqref{shepp_form_joint_density} and \eqref{Corrected_one_Step},  the transition density from $x=s_1$ to $z=s_2$ under the condition  $S(t)<h$ for all $t\in [0,2]$ corrected for discrete time (the corrected form of \eqref{q_x_z}) is given by
\be\label{q_x_z_corrected}
{q}_{h,L}^{(1,2)}(x\to z )& =&
   \frac{1}{{p}_{h,L}^{(0,1)}(x)}\int_{-\infty}^h
   \det  \left(
                                                      \begin{array}{ccc}
                                                       \varphi(s_0)  & \varphi(s_0\!-\!{h_L}\!+\!x) &\varphi(s_0\!-\!2{h_L}\!+\!x\!+\!z) \\
                                                        \varphi({h_L}) & \varphi(x)& \varphi(x\!+\!z\!-\!{h_L}) \\
                                                        \varphi(2{h_L}\!-\!x) & \varphi({h_L})& \varphi(z)
                                                      \end{array}
                                                    \right) ds_0. \nonumber \\
                                                   & =&
   \frac{1}{{p}_{h,L}^{(0,1)}(x)}
   \det  \left(
                                                      \begin{array}{ccc}
                                                       \Phi(h)  & \Phi(h\!-\!{h_L}\!+\!x) &\Phi(h\!-\!2{h_L}\!+\!x\!+\!z) \\
                                                        \varphi({h_L}) & \varphi(x)& \varphi(x\!+\!z\!-\!{h_L}) \\
                                                        \varphi(2{h_L}\!-\!x) & \varphi({h_L})& \varphi(z)
                                                      \end{array}
                                                    \right)
 \ee
 Unlike the transition density \eqref{shepp_form9_correct} (and \eqref{Corrected_one_Step} in the particular case $m=1$), which only depends on $h_L$ and not on $h$, the transition density
 ${q}_{h,L}^{(1,2)}(x\to z )$ depends on both $h$ and $h_L$ and hence the notation. The dependence on $h$ has appeared from integration over the  $s_0 \in (-\infty,h)$.

\subsubsection{Approximations for the BCP ${ P}_{L}(T,h)$}
\label{sec:two_appr}

With  discrete-time corrected transition densities ${q}_h^{(0,1)}(x\to z )$ and ${q}_h^{(1,2)}(x\to z )$, we obtain the corrected versions of the approximations  \eqref{appr_main_1}.\\

\noindent {\bf Approximation 4:}  {\it
$
{ P}_{L}(T,h) \simeq 1-{ F}_{}(1,h,h_L) \cdot \left[{\lambda}_{L,1}(h)\right]^{T-1},
$
 where ${ F}_{}(1,h,h_L)$ is given in \eqref{F_Corrected_shepp_explicit} and ${\lambda}_{L,1}(h)$
is the maximal eigenvalue of the integral operator with kernel $K(x,z)=q^{(0,1)}_{h_L}(x\to z )$ defined in \eqref{Corrected_one_Step}.} \\

\noindent {\bf Approximation 5:}
 {\it
$
{ P}_{L}(T,h) \simeq 1-{ F}_{}(2,h,h_L) \cdot \left[{\lambda}_{L,2}(h)\right]^{T-2},
$
 where ${ F}_{}(2,h,h_L)$ is given in \eqref{F_Corrected_shepp_explicit_2} and ${\lambda}_{L,2}(h)$
is the maximal eigenvalue of the integral operator with kernel $K(x,z)=q^{(1,2)}_{h,L}(x\to z )$ defined in \eqref{q_x_z_corrected}.}\\

Similarly to  $\lambda_m(h)$ from  \eqref{2.27_1},
the maximum eigenvalues ${\lambda}_{L,1}(h)$ and ${\lambda}_{L,2}(h)$ of the operators with kernels $K(x,z)=q^{(0,1)}_{h_L}(x\to z )$ and $K(x,z)=q^{(1,2)}_{h,L}(x\to z )$   are simple and positive; the corresponding eigenfunctions $p(x)$ can be chosen as   probability densities. Both approximations can be used for any $T>0$.

In numerical examples below we approximate  the eigenvalues ${\lambda}_{L,k}(h)$ ($k=1,2$)  using  the methodology described in \cite{Quadrature}, p.154. This methodology is based on the Gauss-Legendre  discretization of  the interval $[-c,h]$, with some large $c>0$, into an $N$-point set $x_1, \ldots, x_N$ (the $x_i$'s are the roots of the \mbox{$N$-th} Legendre polynomial on $[-c,h]$), and the use of the Gauss-Legendre weights $w_i$ associated with points $x_i$; $\lambda_{L,k}(h)$ and $p(x)$ are  then approximated by the largest eigenvalue and associated eigenvector of the matrix
$
D^{1/2}AD^{1/2},
$
where $D = \text{diag}({w}_i)$ and  $A_{i,j} = K(x_i,x_j)$ with the respective kernel $K(x,z)$. If $N$ is large enough then the resulting approximation  to $\lambda_{L,k}(h)$ is arbitrarily accurate. With modern software, computing Approximations 4 and 5 (as well as  Approximation 3) with  high accuracy takes only milliseconds on a regular laptop.

As discussed in the next section, Approximation 5 is more accurate than Approximation 4, especially for small $h$; the accuracies of Approximations 3 and 5 are very similar. Note also that a version of Approximation 4 has been developed in our previous work  \cite{AandZ2019}; this version was based on a different discrete-time approximation (discussed in Section~\ref{T1T2}) of the continuous-time  BCP probability
$P(T,h)$.

\section{Simulation study}\label{sim_study2}

\subsection{Accuracy of approximations for the BCP  ${ P}_{L}(T,h)$  }\label{sim_study_sub}

\begin{figure}[h]
\begin{center}
 \includegraphics[width=0.5\textwidth]{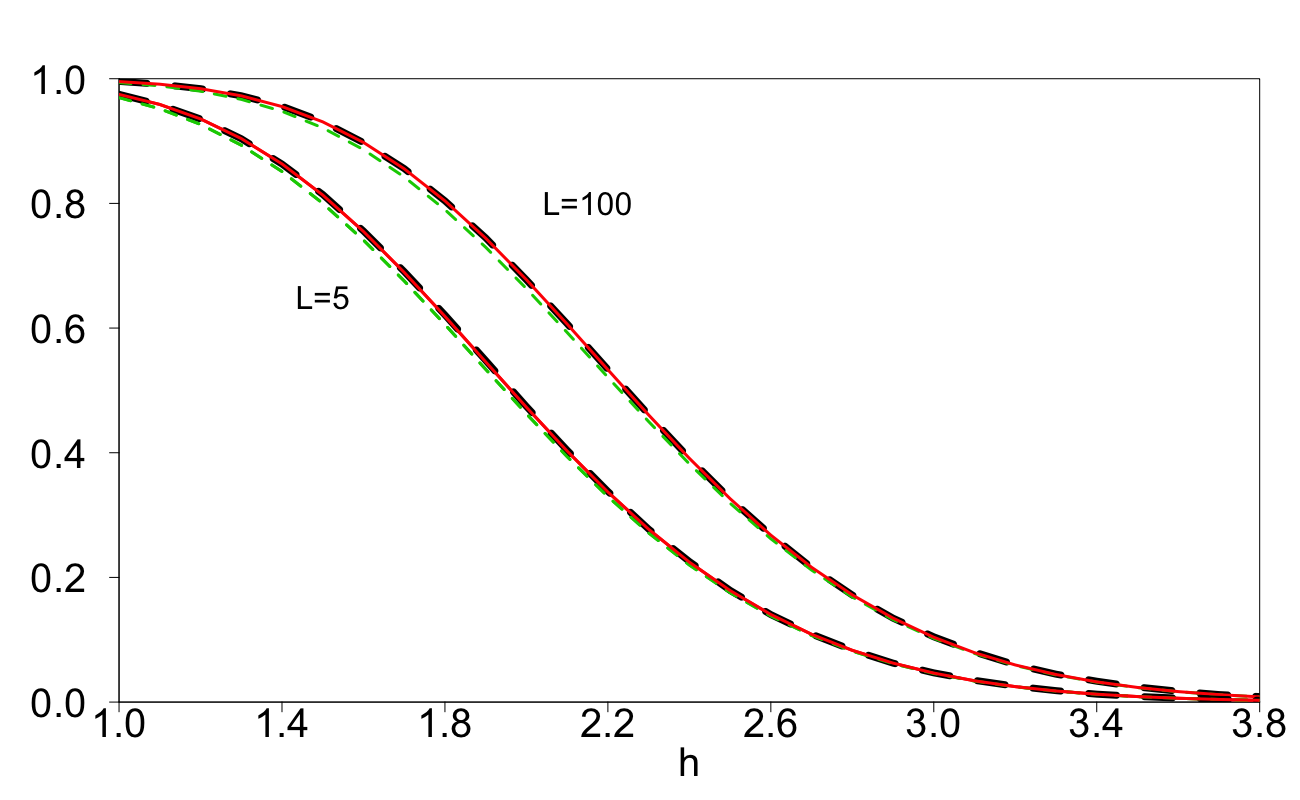}\includegraphics[width=0.5\textwidth]{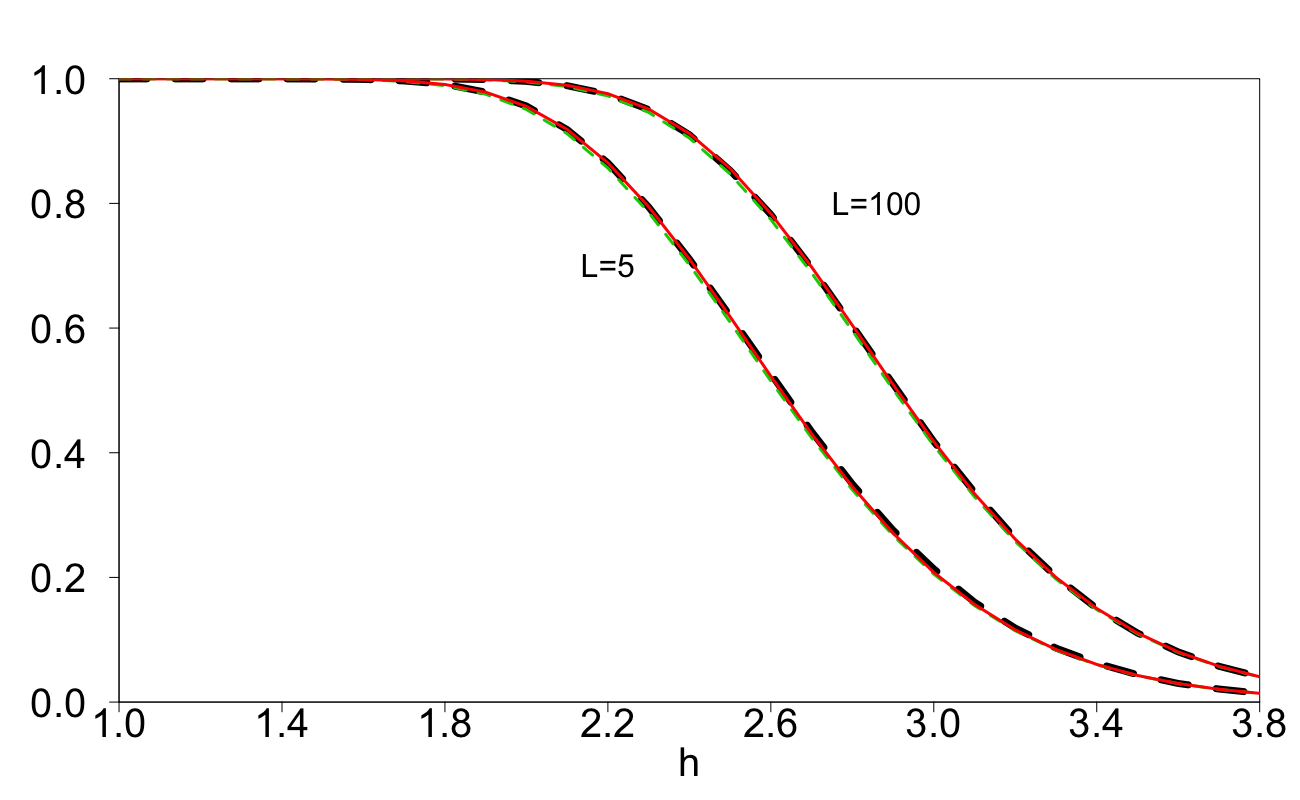}
\end{center}
\caption{Empirical probabilities of reaching the barrier $h$ (dashed black),  Approximation~4 (dashed green) and Approximation 5 (solid red).
Left: $T=10$ with (a) $L=5$ and (b) $L=100$. Right: $T=50$ with (a) $L=5$ and (b)  $L=100$. }
\label{L_10_GSSA}
\end{figure}

In this section we study the quality of Approximations~4 and 5 for the BCP
${ P}_{L}(T,h)$ defined in  \eqref{eq:prob-xi}. Approximation~3 is visually indistinguishable from Approximation~5 and is therefore not plotted (see Table~\ref{lambda_approx}). Without loss of generality, $\varepsilon_j$ in \eqref{eq:sumsq2} are normal r.v.'s with mean $0$ and variance $1$.  The style of Fig.~\ref{L_10_GSSA} is exactly the same as of Fig.~\ref{L_5_CSA}  and is described in
the beginning of Section~\ref{sim_study1}. In  Fig.~\ref{L_10_GSSA}, the dashed green line corresponds to Approximation~4 and the solid red line corresponds to Approximation~5.

From Figure~\ref{L_10_GSSA} we see that the performance of Approximations~4 and 5 is very strong even for small $L$. For small $h$, Approximation~5 is more precise  than Approximation~4  in view of its better accommodation to the non-Markovian nature of the process $S(t)$.

\begin{table}[h]
\centering
\begin{footnotesize}
\begin{tabular}{|c||c|c|c|c|c|c|c|c|c|c|}
\hline
\!\! & $h$=0& $h$=0.5& $h$=1  & $h$=1.5& $h$=2& $h$=2.5& $h$=3& $h$=3.5& $h$=4\\
\hline
\!\!\!\!  ${\lambda}_{L,1}(h)$ \!\!\!\!   &0.28494     &    0.46443  &  0.65331    &   0.81186    & 0.91687     &  0.97090 &    0.99209    & 0.99835   & 0.99974    \\
\!\!\!\! ${\lambda}_{L,2}(h)$ \!\!\!\!  & 0.25744   & 0.43811   & 0.63472    &  0.80239  &  0.91348    & 0.97005   &0.99195     &  0.99833&  0.99974   \\
\!\!\!\!  $\mu_L(h)$ \!\!\!\!   &  0.25527   &  0.43677    &   0.63432 &  0.80241   & 0.91353    &  0.97007  &  0.99195  & 0.99833   &  0.99974    \\
\hline
\end{tabular}
\caption{Values of ${\lambda}_{L,1}(h)$, ${\lambda}_{L,2}(h)$ and ${\mu}_L(h)$  with $L=20$ for different $h$.}
\label{lambda_approx}
\end{footnotesize}
\end{table}

In Table~\ref{lambda_approx}, we display the values of ${\lambda}_{L,1}(h)$, ${\lambda}_{L,2}(h)$ and ${\mu}_L(h)$  with $L=20$ for a number of different $h$. From this table, we see only a small difference between ${\lambda}_{L,2}(h)$ and ${\mu}_L(h)$; this difference is too small to visually differentiate between Approximations~3 and 5 in Fig.~\ref{L_10_GSSA}.

In Tables~\ref{Glaz_error1}, \ref{Glaz_error2} and \ref{Glaz_error3} we numerically compare the performance of Approximations~1 and 3  for approximating ${ P}_{L}(T,h)$ across different values of $L$ and $h$. Since Approximation~1 relies on Monte-Carlo methods, we present the average over 100 evaluations and denote this by $\bar{x}$. We have also provided values for the standard deviation and maximum and minimum of the 100 runs to illustrate the randomised nature of this approximation. These are denoted by  $s$, $Max(x_i) $ and $Min(x_i)$ respectively. The values of ${ P}_{L}(T,h)$ presented in the  tables below are the empirical probabilities of reaching the barrier $h$ obtained by $10^6$ simulations. We have not included Approximation~5 in these tables as results are identical to Approximation~3 up to four decimal places.

\begin{table}[h]
\centering
\begin{footnotesize}
\begin{tabular}{|c||c|c|c|c|c|c|c|c|c|}
\hline
\!\! & $h$=2.5& $h$=2.75& $h$=3  & $h$=3.25& $h$=3.5& $h$=3.75& $h$=4\\
\hline
\!\!\!\! $ \bar{x}$ \!\!  &0.855957 & 0.627299   &   0.376337     &  0.191122    &     0.086253    &   0.033769   &  0.013156    \\
\!\!\!\!  $s$ \!\!   &     0.004127   &    0.008588   &   0.013805   &   0.015181    &   0.012826   &   0.008510  &    0.005131      \\
\!\!\!\!  $Max(x_i) - \bar{x}$ \!\!   &  0.010665     &   0.023748   & 0.029819     &     0.027066  &   0.025629   &  0.016208    &  0.011609       \\
\!\!\!\!  $\bar{x}- Min(x_i)$ \!\!  &  0.012176    & 0.021268   &    0.033211  &  0.041322  &     0.041350  & 0.022650  &    0.018146     \\
\!\!\!\!  Approximation 3\!\!  & 0.854844    &   0.625113  &   0.373863    &   0.188933    &     0.083981       &  0.033833     &    0.012551    \\
\!\!\!\!  ${ P}_{L}(T,h)$  \!\!  & 0.855429    &   0.627463     &     0.376681      &    0.191625        &    0.085697    &   0.034675     &     0.013116      \\
\hline
\end{tabular}
\caption{Average values from 100 evaluations of Approximation 1 for different $h$ along with maximum and minimum with $L=5$ and $T=100$.}
\label{Glaz_error1}
\end{footnotesize}
\end{table}

\begin{table}[h]
\centering
\begin{footnotesize}
\begin{tabular}{|c||c|c|c|c|c|c|c|c|c|}
\hline
\!\! & $h$=2.5& $h$=2.75& $h$=3  & $h$=3.25& $h$=3.5& $h$=3.75& $h$=4\\
\hline
\!\!\!\! $ \bar{x}$ \!\!  & 0.952007 & 0.802073   &   0.554613     &  0.315085     &    0.155331    & 0.066113     &  0.025608    \\
\!\!\!\!  $s$ \!\!   &   0.001479    &  0.004856    &   0.012540   &  0.015050     &   0.015160   & 0.011647    &   0.008129       \\
\!\!\!\!  $Max(x_i) - \bar{x}$ \!\!   & 0.004746      &  0.013360    &  0.027078    &  0.030940        & 0.033991     & 0.024111     &         0.030014\\
\!\!\!\!  $\bar{x}- Min(x_i)$ \!\!  &  0.003662    &  0.010894  &  0.031463    & 0.037715   &  0.041021     & 0.043283  & 0.016997        \\
\!\!\!\!  Approximation 3\!\!  &  0.952475  & 0.802100  &  0.555109   &   0.316076    &    0.153803    & 0.066438    &      0.026143\\
\!\!\!\!  ${ P}_{L}(T,h)$  \!\!  &  0.952818   &   0.803078    &    0.555530      &   0.315784       &     0.153446   &   0.066642    &  0.026244        \\
\hline
\end{tabular}
\caption{Average values from 100 evaluations of Approximation 1 for different $h$ along with maximum and minimum with $L=20$ and $T=100$.}
\label{Glaz_error2}
\end{footnotesize}
\end{table}

\vspace{-0.5cm}
\begin{table}[h]
\centering
\begin{footnotesize}
\begin{tabular}{|c||c|c|c|c|c|c|c|c|c|}
\hline
\!\! & $h$=2.5& $h$=2.75& $h$=3  & $h$=3.25& $h$=3.5& $h$=3.75& $h$=4\\
\hline
\!\!\!\! $ \bar{x}$ \!\!  & 0.979027  & 0.878031  &    0.661247   &   0.402887    & 0.211894    & 0.093329    &  0.039110    \\
\!\!\!\!  $s$ \!\!   &   0.000884    &   0.005502  &  0.014418   &      0.021283  &   0.018493    &  0.020459  &   0.015536     \\
\!\!\!\!  $Max(x_i) - \bar{x}$ \!\!   &   0.001995     &  0.009243   &   0.039695    &     0.040615   &   0.063578    &  0.064306   &    0.037958   \\
\!\!\!\!  $\bar{x}- Min(x_i)$ \!\!  & 0.002414    &  0.020613  &   0.025530 &  0.093876  &    0.038484    &    0.05694   &   0.033748     \\
\!\!\!\!  Approximation 3\!\!  &0.979119   &  0.878481 &   0.660662    &   0.405674 &    0.209313  &   0.094517     &     0.038529   \\
\hline
\end{tabular}
\caption{Average values from 100 evaluations of Approximation 1 for different $h$ along with maximum and minimum with $L=100$ and $T=100$.}
\label{Glaz_error3}
\end{footnotesize}
\end{table}

From Tables~\ref{Glaz_error1}, \ref{Glaz_error2} and \ref{Glaz_error3} we see that with this choice of $T=100$, the errors of approximating $F_L(2,h)$ and $F_L(1,h)$ via the 'GenzBretz' algorithm can accumulate and lead to a fairly significant variation of Approximation~1. This demonstrates the need to average the outcomes of Approximation~1 over a significant number of runs, should one desire an accurate approximation. This may require rather high  computational cost and run time, especially if $L$ is large. On the other hand, evaluation of Approximation~3 is practically instantaneous for all $L$. Even for a very small choice of $L=5$, Table~\ref{Glaz_error1} shows that Approximation~3 still remains very accurate. As $L$ increases from 5 to  $20$, Table~\ref{Glaz_error2} shows that the accuracy of Approximation~3 increases. The averaged Approximation~1 is also very accurate but a larger $L$ appears to produce a larger range for $Max(x_i)$ and $Min(x_i)$ when $h$ is large; this is seen in Table~\ref{Glaz_error3}. Note we have not included empirical values of ${ P}_{L}(T,h)$ in Table~\ref{Glaz_error3} due to the large computational cost.

\subsection{Approximation for the BCP in the case  of non-normal moving sums }
\begin{figure}[!h]
\begin{center}
 \includegraphics[width=0.5\textwidth]{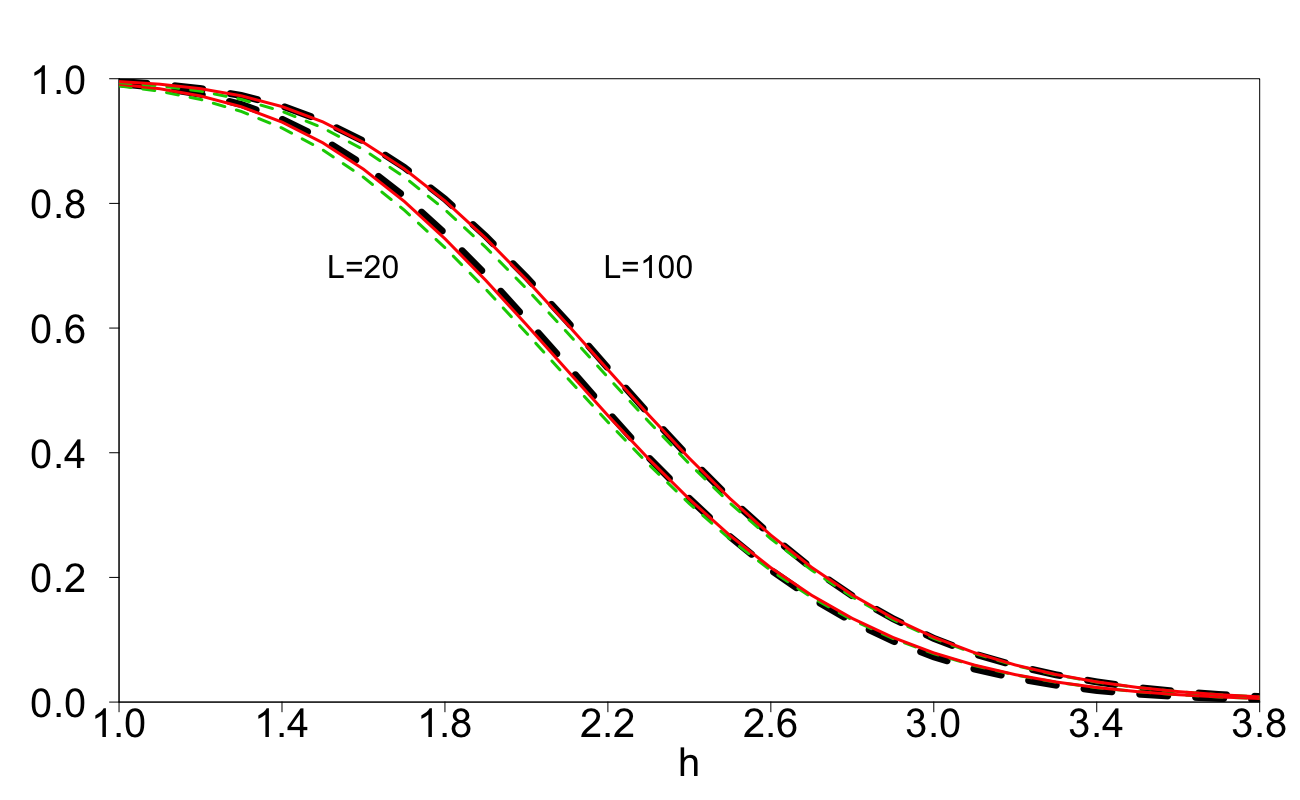}\includegraphics[width=0.5\textwidth]{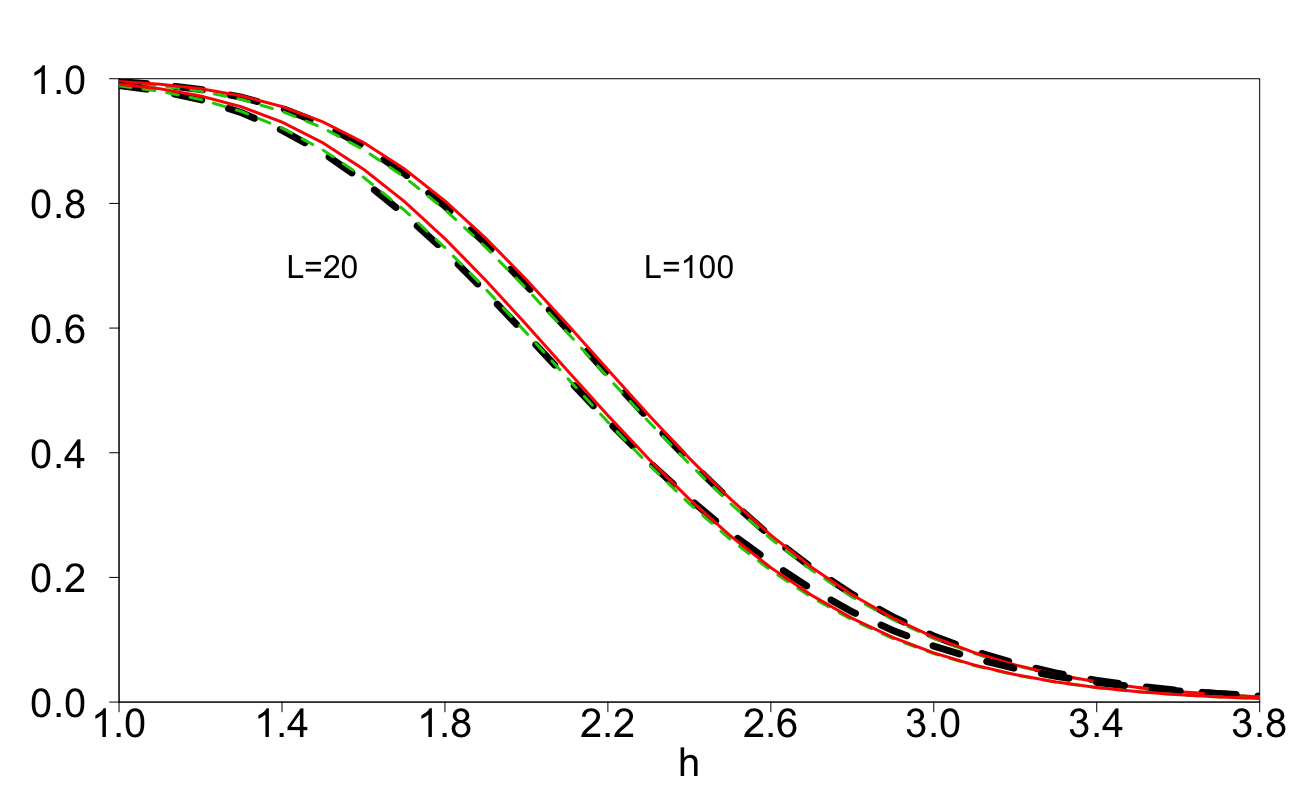}
\end{center}
\caption{Empirical probabilities of reaching the barrier $h$ (dashed black),  Approximation~4 (dashed green) and Approximation 5 (solid red).
Left: $\varepsilon_i\sim \text{Uniform}[0,1]$ and $T=10$ with (a) $L=20$ and (b) $L=100$. Right: $\varepsilon_i\sim \text{Laplace}[0,1]$ and $T=10$ with (a) $L=20$ and (b)  $L=100$. }
\label{fig:non_normal}
\end{figure}

 Approximations~3, 4 and 5 remain very accurate when then the original $\varepsilon_i$ in \eqref{eq:sumsq2} are not exactly normal. We consider two cases: (a)  $\varepsilon_i$ are uniform r.v's on [0,1] and (b)  $\varepsilon_i$ are Laplace r.v's with mean zero and scale parameter 1. Simulation results are shown in Figure~\ref{fig:non_normal}; this figure has the same style as  figures in Sections~\ref{CSA_sim_study} and \ref{sim_study_sub}.

Some  selected values used for plots in Figure~\ref{fig:non_normal}  are:\\
   $h=2,L=20\!:\;\;$ Emp: $0.6045 \pm 0.0030\;(0.6123 \pm 0.0030) \; [0.5894 \pm 0.003];\;\,$ Ap. 4(5): 0.5921(0.6054); \\
 $h=2,L=\!100\!:\;$ Emp: $0.6771 \pm 0.0029\;(0.6801 \pm 0.0029)\;[0.6722 \pm 0.003 ];\;\,$ Ap. 4(5): 0.6633(0.6775);\\
   $h=3,L=20\!:\;\;$ Emp: $0.0788 \pm 0.0017\;(0.0710\pm  0.0016)\; [0.0915 \pm 0.002];\;\,$ Ap. 4(5): 0.0777(0.0789);\\
   $h=3,L=\!100\!:\;$ Emp: $0.1039 \pm 0.0019\;(0.1033 \pm 0.0019)\; [0.1048 \pm 0.002];\;\,$ Ap. 4(5): 0.1022(0.1034).

Here we  provided  means and 95\% confidence intervals for the empirical (Emp) values of the BCP ${ P}_{L}(T,h)$ (with $T=M/L=10$) computed from  100\,000 Monte-Carlo runs of the sequences of the moving sums \eqref{eq:sumsq2} with normal (no brackets), uniform (regular brackets) and Laplace (square brackets) distributions for  $\varepsilon_i$ in \eqref{eq:sumsq2}. Values of Approximations (Ap.) 4 and 5 are also given.

From Figure~\ref{fig:non_normal} and associated numbers we can make the following conclusions: (a) the BCP ${ P}_{L}(T,h)$ for the case where $\varepsilon_i$ in \eqref{eq:sumsq2} are uniform is closer to the case where $\varepsilon_i$ are normal, than for the case where $\varepsilon_i$ have  Laplace distribution; (b) as $L$ increases, the probabilities ${ P}_{L}(T,h)$  in the cases of uniform and Laplace distributions of $\varepsilon_i$ become  closer to the BCP for the case of normal $\varepsilon_i$ and hence the  approximations  to the BCP become more precise; (c) accuracy of Approximation~5 is excellent for the case of normal $\varepsilon_i$ and remains  very good in the case of uniform $\varepsilon_i$; it is also rather good in the case when $\varepsilon_i$ have Laplace distribution; (d) Approximation~4 is slightly less accurate than Approximation 5 (and Approximation 3) for the case of normal and uniform $\varepsilon_i$ (this is in a full agreement with discussions in Sections~\ref{sec:two_appr} and \ref{sim_study_sub}); however, Approximation~4 is very simple and can still be considered as rather accurate.

\subsection{Approximation for the BCP in the case  of moving weighted sums }

We have also investigated the performance of Approximation~5  (and 3) after introducing particular weights into  \eqref{eq:sumsq2}.
  We explored the following two ways of incorporating weights:
\begin{itemize}
  \item[(i)]
  $L$ random weights $w_1,w_2,\ldots,w_L$, with $w_i$ i.i.d. uniform on $[0, 2]$, are associated with a position in the moving window; this results in the moving weighted sum
\begin{eqnarray*}
S_{n,w,L}:= \sum_{j=n+1}^{n+L}w_{j-n} \varepsilon_j\, \;\; (n=0,1, \ldots,M)\, ;
\end{eqnarray*}
  \item[(ii)]
 $M+L$ random weights $w_1,\ldots,w_{M+L}$ are associated with r.v. $\varepsilon_1,\ldots,\varepsilon_{M+L}$; here  $w_j$ are i.i.d. uniform r.v's on [0,2]; this gives the moving weighted sum
\begin{eqnarray*}
S_{n,w,L}:= \sum_{j=n+1}^{n+L}w_{j} \varepsilon_j\, \;\; (n=0,1, \ldots,M).
\end{eqnarray*}

  \end{itemize}

\begin{figure}[h]
\begin{center}
 \includegraphics[width=0.5\textwidth]{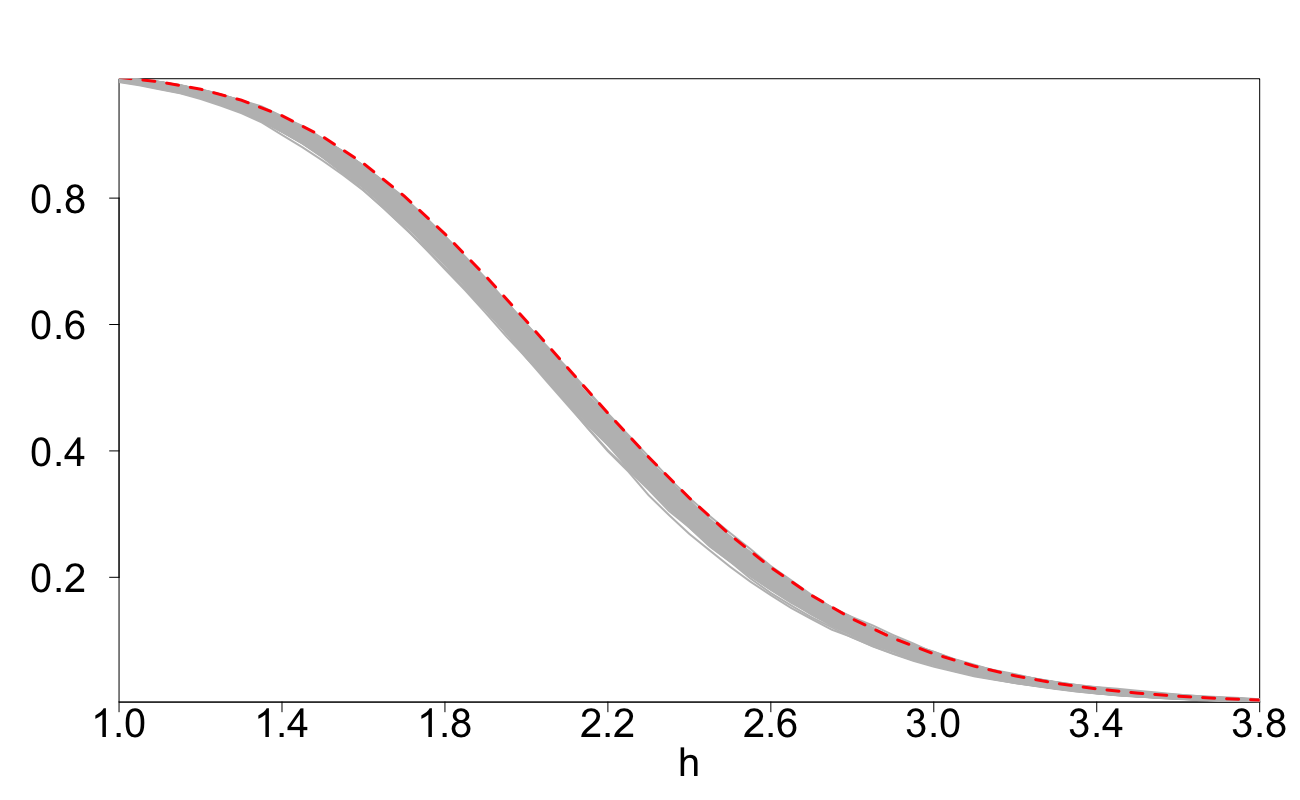}\includegraphics[width=0.5\textwidth]{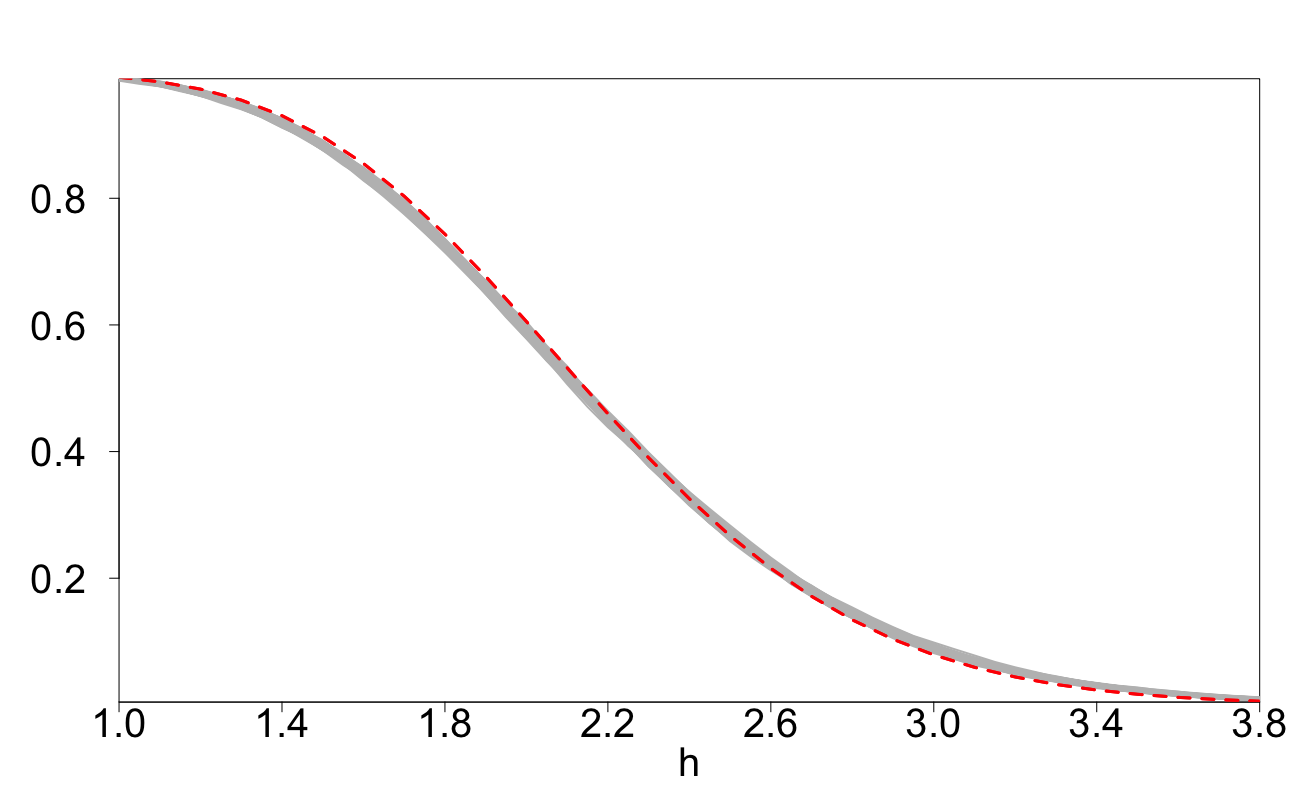}
\end{center}
\caption{BCP for the weighted sums (grey) against Approximation~5 for the BCP for non-weighted moving sums (red dotted line).  Left: case (i) with $L=20,M=200, T=10$. Right: case (ii)  with $L=20,M=200,T=10$.}
\label{fig:moving_weights}
\end{figure}

Simulations results are shown in Fig.~\ref{fig:moving_weights}.
In both cases, we have repeated simulations 1,000 times and plotted all the curves representing the BCP as functions of $h$  in grey colour and Approximation~5 for the BCP for the non-weighted case (when all weights $w_j=1$)  as red dashed line. We can see that for both scenarios the Approximation~5 for the BCP in the non-weighted case gives fairly accurate approximation for the weighted BCP. Similar results have been observed for other values of $L$ and $T$.

\section{Approximating Average Run Length (ARL) }
\label{ARL_section}

In this section, we provide approximations to the probability distribution of  the moment of time
  \mbox{$
  \tau_{H}(\s)\!:=\!\min \{ n\!\geq \! 0\!: S_{n,L} \!\geq \! H \}
  $}
   when the sequence $\s=\{S_{0,L},S_{1,L}, \ldots\}$  reaches  the threshold $H$ for the first time. Note that $\tau_{H}(\s)=\tau_{h}(\mathbb{X})$, where   \mbox{$
  \tau_{h}(\mathbb{X})\!:=\!\min \{ n\!\geq \! 0\!: \xi_{n,L} \!\geq \! h \}
  $} and $\mathbb{X}=\{\xi_{0,L},\xi_{1,L}, \ldots\}$.
  The BCP ${\cal P}_{\s}(M,H,L)$, considered as a function of $M$, is the c.d.f. of this probability distribution: ${\cal P}_{\s}(M,H,L)=
  {\rm Pr}\left(  \tau_{H}(\s) \leq M \right)$.
The average run length (ARL) until $\s$ reaches $H$ for the first time is
\be
\label{arl}
{\rm ARL}_H(\s):= \sum_{n=0}^\infty n {\rm Pr} \{ \tau_H=n\}= \int_0^{\infty}{M d {\cal P}_{\s}(M,H,L) }  \, .
\ee

Note that ${\rm ARL}_H(\s)={\rm ARL}_h(\mathbb{X})$.
The diffusion approximation to the time moment  $\tau_{h}(\mathbb{X})$ is   $\tau_{h}(S(t)):=\min \{ {t\geq 0}:\;  S(t) \geq h \}
  $, which is
   the time moment when the process $S(t)$ reaches $h$. The distribution of $\tau_{h}(S(t))$
  has the form:
  \bea
  (1-\Phi(h)) \delta_0 (d s) + q(s,h,S(t)) d s \, ,s \geq 0,
  \eea
   where $\delta_0 (d s)$ is the delta-measure concentrated at 0 and
  \begin{equation}
\label{eq:first_pass_time}
q(s,h,S(t))=\frac{d}{ds}{P}_{}(s,h),\;\;\;0<s<\infty\, .
\end{equation}
The function $q(s,h,S(t))/\Phi(h)$, considered as a function of $s$,  is a probability density function on $(0, \infty)$ since
\bea
\int_0^\infty q(s,h,S(t)) ds = 1- {P}_{}(0,h)= \Phi(h)\, .
\eea
From this, the diffusion approximation for ${\rm ARL}_H(\mathbb{X})/L$ is
\begin{equation}
{\rm ARL}_h(S(t)) =  \E (\tau_{h}(S(t)))=\int_0^{\infty}{s \,q(s,h,S(t))ds}\, .\label{ARL_form}
\end{equation}

The diffusion approximation  \eqref{ARL_form} should be corrected for discrete time; otherwise it is poor, especially for small $L$.
As shown in Section~\ref{sim_study2}, Approximations~3 and 5 are very accurate approximations for ${ P}_{L}(T,h)$ and can be used
for all  $T >0$.
 We shall use Approximation~3 to formulate our approximations but note that the use of Approximation~5 would give very similar results.

We define the  approximation $\hat{q}(s,h)$ for the probability density function of $\tau_{h}(\mathbb{X})$$/L$ by
\begin{eqnarray*}\label{q_t_h}
\hat{q}(s,h)&=&
\frac{d}{ds} \left\{ 1- { F}_{}(2,{h},h_L) \cdot  \mu_L(h)^{s-2}  \right\}
  = - { F}_{}(2,{h},h_L)\log\left (\mu_L(h) \right )\cdot    \mu_L(h) ^{s-2}, \,\,\,\, s >0.
\end{eqnarray*}
The corresponding approximation  for ${\rm ARL}_h(\mathbb{X})$  is
\begin{eqnarray}
{\rm ARL}_h(\mathbb{X})  =  \E \tau_{h}(\mathbb{X}) &\cong& L \int_0^{\infty}{s\hat{q}(s,h) ds}
 \,= -\frac{L\cdot   { F}_{}(2,{h},h_L)}{\mu_L(h)^2\log( \mu_L(h))}.\label{ARL_form_app}
\end{eqnarray}

The standard deviation of $\tau_{h}(\mathbb{X})$, denoted $SD(\tau_{h}(\mathbb{X}))$, is approximated by:
\begin{equation}\label{SD_formula}
SD(\tau_{h}(\mathbb{X})) \cong L \,\left[\int_{0}^{\infty}s^2 \, \hat{q}(s,h) ds - \left(\int_{0}^{\infty}s\,\hat{q}(s,h) ds \right)^2 \right]^{1/2}.
\end{equation}

In this paper, we define ARL in terms of the number of random variables $\xi_{n,L}$ rather than number of random variables $\varepsilon_j$. This means we have to modify the   approximation for ARL of \cite{Glaz2012} by subtracting $L$. The standard deviation approximation in \cite{Glaz2012} is not altered.

%
%

The Glaz approximations for ${\rm ARL}_h(\mathbb{X})$ and $SD(\tau_{h}(\mathbb{X}))$ are as follows:
\be \label{Glaz_ARL}
 \E _G(\tau_{h}(\mathbb{X})) = \! \sum_{j=L}^{2L}{ F}_{L}((j/L-1),h) \!+\! \frac{{F}_{L}(1,h)}{{F}_{L}(1,h)-{F}_{L}(2,h)}\sum_{j=1}^{L}({ F}_{L}(1+j/L,h))\, , \;\;\;
\ee
\be \label{Glaz_SD}
SD_G(\tau_{h}(\mathbb{X})) = \bigg [&&\!\!\!\!\!\!\! L(L-1) \! + \!2\sum_{j=L}^{3L}j({ F}_{L}(j/L-1,h)) + \frac{2Lx(3-2x)}{(1-x)^2}\sum_{j=1}^{L}{ F}_{L}(1+j/L,h) \nonumber \\
&+& \frac{2x}{1-x}\sum_{j=1}^{L}j({ F}_{L}(1+j/L,h)) \!+\!  \E _G(\tau_{h}(\mathbb{X})) \!-\!  \E _G(\tau_{h}(\mathbb{X}))^2 \bigg]^{1/2}\, , \;\;\;\;\;\;\;\;\;
\ee
where $x={ F}_{L}(2,h) /{F}_{L}(1,h)$.

 In Tables~\ref{expected_run_length} and \ref{SD_run_length} we assess the accuracy of  the approximations  \eqref{ARL_form_app} and \eqref{SD_formula}
 and also Glaz approximations \eqref{Glaz_ARL} and \eqref{Glaz_SD}. In these tables, the values of ${\rm ARL}_h(\mathbb{X})$ and $SD(\tau_{h}(\mathbb{X}))$ have been calculated using $100,000$ simulations. Since the Glaz approximations rely on Monte Carlo methods,
 in the tables we have reported value $2s$-confidence intervals computed from 150 evaluations.


\begin{table}[!h]
\centering
\caption{Approximations for ${\rm ARL}_h(\mathbb{X})$ and $SD(\tau_{h}(\mathbb{X}))$ with $L=10$.}
\label{expected_run_length}
\begin{tabular}{|c||c|c|c|c|c|c|c|}
\hline
$h$ & 2 & 2.25 &2.5 &2.75 & 3 &3.25 & 3.5  \\
\hline
 \eqref{ARL_form_app}& 126 &217  &395   & 759  &  1551  & 3375  &  7837           \\
\eqref{Glaz_ARL}  & 126 \!$\pm$ \!\!1   & 218 \!$\pm$ \!\!2 & 394 \!$\pm$ \!\!5  & 756 \!$\pm$ \!\!17  &  1545 \!$\pm$ \!\!65 & 3388 \!$\pm$ \!\!300 &   7791 \!$\pm$ \!\!1100              \\
${\rm ARL}_h(\mathbb{X})$ & 127 & 218 & 396  &757  & 1550  & 3344  &7721      \\
\hline
\end{tabular}

\medskip

\centering

\begin{tabular}{|c||c|c|c|c|c|c|c|}
\hline
$h$ & 2 & 2.25 &2.5 &2.75 & 3 &3.25 & 3.5  \\

\hline
 \eqref{SD_formula}&  129  & 220   &  397   &   761   & 1553    &  3377  &    7839          \\
\eqref{Glaz_SD}  & 129 \!$\pm$ \!\!1  &  220 \!$\pm$ \!\!2  & 397 \!$\pm$ \!\!5  &  758 \!$\pm$ \!\!17  &  1549 \!$\pm$ \!\!65 &     3389 \!$\pm$ \!\!300   &7793 \!$\pm$ \!\!1100  \\
\!\!\! $SD(\tau_{h}(\mathbb{X})) \!\!$ &129   & 221  & 395  &758    & 1550    & 3341     & 7716            \\
\hline
\end{tabular}
\end{table}

\begin{table}[h]
\centering
\caption{Approximations for ${\rm ARL}_h(\mathbb{X})$ and $SD(\tau_{h}(\mathbb{X}))$ with $L=50$.}
\label{SD_run_length}

\begin{tabular}{|c||c|c|c|c|c|c|c|}
\hline
$h$ & 2 & 2.25 &2.5 &2.75 & 3 &3.25 & 3.5   \\
\hline
 \eqref{ARL_form_app}&  471 & 791  & 1392   & 2587  & 5099    &  10695  &    23918           \\
\eqref{Glaz_ARL}  & 471 \!$\pm$ \!\!3   &  791 \!$\pm$ \!\!7   &   1393 \!$\pm$ \!\!25   & 2597 \!$\pm$ \!\!75  & \!\!\! 5101 \!$\pm$ \!\!270 \!\!\! & \!\!\!\! 10708 \!$\pm$ \!\!1250 \!\!\!\! & \!\!\!\!   24639 \!$\pm$ \!\!5800  \!\!\!\!        \\
${\rm ARL}_h(\mathbb{X})$ & 472  &792  & 1397  &  2588 & 5085  & 10749   &    24131      \\
\hline
\end{tabular}
\bigskip

\centering

\begin{tabular}{|c||c|c|c|c|c|c|c|}
\hline
$h$ & 2 & 2.25 &2.5 &2.75 &3 &3.25 & 3.5   \\

\hline
 \eqref{SD_formula}& 485   & 804 &  1404  &   2598  &  5109   & 10704   &       23924         \\
\eqref{Glaz_SD}  & 481 \!$\pm$ \!\!3  &  802 \!$\pm$ \!\!7  &  1404 \!$\pm$ \!\!25  & 2608 \!$\pm$ \!\!75  & \!\!\!\! 5147 \!$\pm$ \!\!270 \!\!\!\! & \!\! 10716 \!$\pm$ \!\!1250 \!\!\! \!\!\! & \!\!\! 24649 \!$\pm$ \!\!5800  \!\!\!\!        \\
\!\!$SD(\tau_{h}(\mathbb{X})) \!\!$&485   &  804 & 1407  & 2600  &  5093 &  10762   &  24105      \\
\hline
\end{tabular}

\end{table}

 Tables~\ref{expected_run_length} and \ref{SD_run_length} show that the approximations developed in this paper perform strongly and are similar, for small or moderate $h$, to the Glaz approximations.  For \mbox{$h\!\geq \! 3$,} the Glaz approximation produces rather large uncertainty intervals and the uncertainty quickly deteriorates with the increase of $h$.  This is due to the fairly large uncertainty intervals formed by Approximation~1 when approximating ${ P}_{L}(T,h)$ with large $h$ and hence small ${ P}_{L}(T,h)$, as discussed in Section~\ref{sim_study_sub}. The approximations developed in this paper are deterministic and are much simpler in comparison to the Glaz approximations. Moreover, they  do not deteriorate for large  $h$.

\section*{Acknowledgment} The authors are grateful to the referees for careful reading of the manuscript and useful comments.

\bibliographystyle{spmpsci}

\bibliography{changepoint1}

\begin{thebibliography}{10}
\providecommand{\url}[1]{{#1}}
\providecommand{\urlprefix}{URL }
\expandafter\ifx\csname urlstyle\endcsname\relax
  \providecommand{\doi}[1]{DOI~\discretionary{}{}{}#1}\else
  \providecommand{\doi}{DOI~\discretionary{}{}{}\begingroup
  \urlstyle{rm}\Url}\fi

\bibitem{Bau2}
Bauer, P., Hackl, P.: An extension of the {MOSUM} technique for quality
  control.
\newblock Technometrics \textbf{22}(1), 1--7 (1980)

\bibitem{Chu}
Chu, C.S.J., Hornik, K., Kaun, C.M.: {MOSUM} tests for parameter constancy.
\newblock Biometrika \textbf{82}(3), 603--617 (1995)

\bibitem{eiauer1978use}
Eiauer, P., Hackl, P.: The use of {MOSUMS} for quality control.
\newblock Technometrics \textbf{20}(4), 431--436 (1978)

\bibitem{genz2009computation}
Genz, A., Bretz, F.: Computation of Multivariate Normal and t Probabilities.
\newblock Lecture Notes in Statistics. Springer-Verlag, Heidelberg (2009)

\bibitem{GenzR}
Genz, A., Bretz, F., Miwa, T., Mi, X., Leisch, F., Scheipl, F., Hothorn, T.:
  {mvtnorm}: Multivariate Normal and t Distributions (2018).
\newblock \urlprefix\url{https://CRAN.R-project.org/package=mvtnorm}.
\newblock R package version 1.0-8:
  \textit{`https://CRAN.R-project.org/package=mvtnorm'}

\bibitem{Glaz_old}
Glaz, J., Johnson, B.: Boundary crossing for moving sums.
\newblock Journal of {A}pplied Probability \textbf{25}(1), 81--88 (1988)

\bibitem{Glaz2012}
Glaz, J., Naus, J., Wang, X.: Approximations and inequalities for moving sums.
\newblock Methodology and Computing in Applied Probability \textbf{14}(3),
  597--616 (2012)

\bibitem{glaz1991tight}
Glaz, J., Naus, J.I.: Tight bounds and approximations for scan statistic
  probabilities for discrete data.
\newblock The Annals of Applied Probability \textbf{1}(2), 306--318 (1991)

\bibitem{glaz2001scan}
Glaz, J., Naus, J.I., Wallenstein, S., Wallenstein, S., Naus, J.I.: Scan
  statistics.
\newblock Springer (2001)

\bibitem{glaz2009scan2}
Glaz, J., Pozdnyakov, V., Wallenstein, S.: Scan Statistics: Methods and
  Applications.
\newblock Birkhäuser, Boston (2009)

\bibitem{Haiman}
Haiman, G.: First passage time for some stationary processes.
\newblock Stochastic Processes and their Applications \textbf{80}(2), 231--248
  (1999)

\bibitem{Quadrature}
Mohamed, J., Delves, L.: Computational Methods for Integral Equations.
\newblock Cambridge University Press (1985)

\bibitem{MZ2003}
Moskvina, V., Zhigljavsky, A.: An algorithm based on {S}ingular {S}pectrum
  {A}nalysis for change-point detection.
\newblock Communications in Statistics---Simulation and Computation
  \textbf{32}(2), 319--352 (2003)

\bibitem{AandZ2019}
Noonan, J., Zhigljavsky, A.: Approximations for the boundary crossing
  probabilities of moving sums of normal random variables.
\newblock Communications in Statistics-Simulation and Computation, 1--22
  (2019)

\bibitem{ReedSimon}
Reed, M., Simon, B.: Methods of Modern Mathematical Physics: Scattering theory
  Vol. 3.
\newblock Academic Press (1979)

\bibitem{Shepp71}
Shepp, L.: First passage time for a particular {G}aussian process.
\newblock The Annals of Mathematical Statistics \textbf{42}(3), 946--951 (1971)

\bibitem{Sieg_book}
Siegmund, D.: Sequential Analysis: Tests and Confidence Intervals.
\newblock Springer Science \& Business Media (1985)

\bibitem{Sieg_paper}
Siegmund, D.: Boundary crossing probabilities and statistical applications.
\newblock The Annals of Statistics \textbf{14}(2), 361--404 (1986)

\bibitem{slepian1961first}
Slepian, D.: First passage time for a particular {G}aussian process.
\newblock The Annals of Mathematical Statistics \textbf{32}(2), 610--612 (1961)

\bibitem{waldmann1986bounds}
Waldmann, K.H.: Bounds to the distribution of the run length in general
  quality-control schemes.
\newblock Statistische Hefte \textbf{27}(1), 37 (1986)

\bibitem{wang2014variable}
Wang, X., Glaz, J.: Variable window scan statistics for normal data.
\newblock Communications in Statistics-Theory and Methods \textbf{43}(10-12),
  2489--2504 (2014)

\bibitem{wang2014multiple}
Wang, X., Zhao, B., Glaz, J.: A multiple window scan statistic for time series
  models.
\newblock Statistics \& Probability Letters \textbf{94}, 196--203 (2014)

\bibitem{Xia}
Xia, Z., Guo, P., Zhao, W.: Monitoring structural changes in generalized linear
  models.
\newblock Communications in Statistics—Theory and Methods \textbf{38}(11),
  1927--1947 (2009)

\end{thebibliography}


\end{document}